\documentclass[11pt,a4paper]{amsart}
\usepackage[T1]{fontenc}
\usepackage[utf8]{inputenc}
\usepackage{lmodern}
\usepackage[margin=28mm]{geometry}
\usepackage{amsmath,amssymb,amsthm,mathtools}
\usepackage{booktabs}
\usepackage{enumitem}
\usepackage{listings}
\usepackage{microtype}
\usepackage{xcolor}
\usepackage[hidelinks]{hyperref}

\lstdefinestyle{mma}{
	language=Mathematica,
	basicstyle=\ttfamily\footnotesize,
	upquote=true, columns=fullflexible, keepspaces=true,
	showstringspaces=false, tabsize=2,
	breaklines=true, breakatwhitespace=true,
	postbreak=\mbox{\textcolor{gray}{$\hookrightarrow$}\space},
	xleftmargin=2em, frame=single, framexleftmargin=1.5em, framesep=4pt,
	rulecolor=\color{gray!50},
	commentstyle=\itshape\color{gray!70!black},  % comments: grey italic
	stringstyle=\color{green!40!black},          % labels: dark green
	keywordstyle={},      
	aboveskip=8pt, belowskip=8pt,
}

\numberwithin{equation}{section}
\newtheorem{theorem}{Theorem}[section]
\newtheorem{proposition}[theorem]{Proposition}
\newtheorem{lemma}[theorem]{Lemma}
\newtheorem{corollary}[theorem]{Corollary}
\theoremstyle{definition}

\theoremstyle{remark}
\newtheorem{remark}[theorem]{Remark}
\newcommand{\R}{\mathbb R}
\newcommand{\Sph}{\mathbb S}
\newcommand{\eps}{\epsilon}
\newcommand{\grad}{\operatorname{grad}}

\newcommand{\tr}{\operatorname{trace}}
\newcommand{\Id}{\operatorname{Id}}
\newcommand{\Res}{\operatorname{Res}}
\newcommand{\Gf}{\mathfrak G}
\newcommand{\Zf}{\mathfrak Z}
\newcommand{\tE}{\widetilde E}

\title[Biharmonic surfaces with flat normal bundle in $\Sph^4$]{Biharmonic surfaces with flat normal bundle\\ in the four-dimensional sphere}
\date{}
\author{\c Stefan Andronic, Stefano Montaldo, Cezar Oniciuc and Antonio Sanna}

\address{Faculty of Mathematics, Al. I. Cuza University of Iasi, Blvd. Carol I, no. 11, 700506 Iasi, Romania} \email{stefanandronic215@gmail.com}

\address{Dipartimento di Matematica e Informatica, Universit\` a degli Studi di Cagliari, Via Ospedale 72, 09124 Cagliari, Italia} \email{montaldo@unica.it}

\address{Faculty of Mathematics, Al. I. Cuza University of Iasi, Blvd. Carol I, no. 11, 700506 Iasi, Romania} \email{oniciucc@uaic.ro}

\address{Dipartimento di Matematica e Informatica, Universit\` a degli Studi di Cagliari, Via Ospedale 72, 09124 Cagliari, Italia}

\email{antonio.sanna4@unica.it}


\thanks{The authors S.M. and A.S. are members of the Italian National Group G.N.S.A.G.A. of INdAM. The author A. S. was supported by a NRRP scholarship - funded by the European Union - NextGenerationEU - Mission 4, Component 1, Investment 3.4. The work was partially supported by the Project {PRoBIKI} funded by Fondazione di Sardegna.
}

\begin{document}
\begin{abstract}
We prove that a proper biharmonic surface with flat normal bundle in the
four-dimensional unit sphere $\Sph^4$ lies minimally in a small hypersphere of
radius $1/\sqrt2$. Since the surfaces with constant mean curvature and those with
parallel normalized mean curvature vector field are already understood, the core of
the paper is the non-existence of biharmonic surfaces with flat normal bundle in
$\Sph^4$ whose mean curvature is not constant and whose normalized mean curvature
vector field is not parallel. We first show, adapting an argument of Fu, Hong and
Yang, that the gradient of the mean curvature is a principal direction. Then the
problem reduces to a first order system of ordinary differential equations with a
constraint, and we show that this system has no admissible solutions by an
elimination argument.
\end{abstract}
\maketitle

\section{Introduction}\label{sec:introduction}

Harmonic maps, which are critical points of the energy functional, provide a
variational framework in which the geometry of maps $\varphi:(M,g)\to(N,h)$ between
Riemannian manifolds and the theory of minimal immersions can be studied together.
When $M$ is compact, the energy of $\varphi$ is
\[
 E(\varphi)=\frac12\int_M|d\varphi|^2\,dv_g .
\]
The Euler--Lagrange equation associated to the energy functional is the vanishing of
the tension field
\[
 \tau(\varphi)=\tr_g\nabla d\varphi,
\]
that is, the trace of the covariant derivative of the differential of the map
$\varphi$. In particular, the tension field of an isometric immersion of an
$m$-dimensional manifold is $mH$, where $H$ is the mean curvature vector field.
Thus isometric immersions are harmonic precisely when they are minimal. Biharmonic
maps arise by replacing the energy with the bienergy
\[
 E_2(\varphi)=\frac12\int_M|\tau(\varphi)|^2\,dv_g .
\]
In \cite{Jiang} Jiang established the first and second variation formulas for this
functional. The resulting Euler--Lagrange equation is
\[
 \tau_2(\varphi)=-\Delta^\varphi\tau(\varphi)
 -\tr_g R^N\bigl(d\varphi,\tau(\varphi)\bigr)d\varphi=0,
\]
a fourth-order semi-linear elliptic system. Every harmonic map is biharmonic, and a
biharmonic map which is not harmonic is called \emph{proper biharmonic}. We refer to
\cite{FetcuOniciuc2022, ChenOu2020Book} for an account of the basic theory and its relation to submanifold geometry.

When the target is the Euclidean space, an isometric immersion is biharmonic if and
only if its position vector satisfies $\Delta^2\varphi=0$, and Chen's conjecture
predicts that every biharmonic submanifold of the Euclidean space is minimal
\cite{Chen}. Very recently, Fu, Hong and Yang proved Chen's conjecture for surfaces in
$\R^4$ with flat normal bundle \cite{FHY}. The situation is different when the ambient space is the unit Euclidean
sphere $\Sph^n$, which carries many proper biharmonic submanifolds. The basic example
is the small hypersphere $\Sph^{n-1}(2)$ of radius $1/\sqrt2$, where, as in \cite{AMOS},
$\Sph^m(c)$ denotes the $m$-dimensional sphere of constant sectional curvature $c$, and,
more generally, a minimal submanifold of $\Sph^{n-1}(2)$ is a proper biharmonic
submanifold of $\Sph^n$ (see \cite{CMO}). The known examples of proper biharmonic
submanifolds of spheres have constant mean curvature, and it was conjectured in
\cite{BMO} that every proper biharmonic submanifold of $\Sph^n$ has constant mean
curvature. For surfaces in $\Sph^4$ with
constant mean curvature (CMC) the picture is complete: Balmu\c{s} and Oniciuc proved
in \cite{BO} that a proper biharmonic CMC surface in $\Sph^4$ lies minimally in a
small hypersphere $\Sph^3(2)$, of radius $1/\sqrt2$. In
\cite[Conjecture~2]{AMOS} we conjectured that every proper biharmonic surface in
$\Sph^4$ lies minimally in a small hypersphere $\Sph^3(2)$, that is, that the CMC
hypothesis in \cite{BO} can be dropped.

For an isometric immersion $\varphi:M^2\to\Sph^4$, the bitension field has a normal
and a tangent component. With the sign convention $\Delta=-\tr\nabla^2$ and a local
orthonormal frame $\{e_1,e_2\}$ of $TM$, the immersion $\varphi$ is biharmonic if and
only if
\[
 \Delta^\perp H+\sum_{i=1}^2B(e_i,A_He_i)-2H=0
 \qquad\text{and}\qquad
 2\sum_{i=1}^2A_{\nabla^\perp_{e_i}H}e_i+\grad|H|^2=0,
\]
where $B$ is the second fundamental form, $A$ the shape operator and $\nabla^\perp$,
$\Delta^\perp$ the connection and the Laplacian of the normal bundle. The tangent
equation has an independent variational meaning: it is equivalent to the vanishing
of the divergence of the stress-bienergy tensor \cite{CaddeoMontaldoOniciucPiu2014, Jiang2, LMO}, and the immersions
which satisfy it are called \emph{biconservative}. On the set where $H\neq0$, we write
$H=fE_3$, where $f=|H|$ is the mean curvature function and $E_3$ is the
\emph{normalized mean curvature vector field}. Then
$\nabla^\perp_XH=X(f)E_3+f\nabla^\perp_XE_3$. For a hypersurface the second term
vanishes, while in codimension two it does not, and this is where the normal
connection enters the biharmonic equation. When $\nabla^\perp E_3=0$, the surface is
said to have \emph{parallel normalized mean curvature vector field} (PNMC). In
codimension two, a PNMC surface has flat normal bundle, since $E_3$ and the unit normal
orthogonal to it are both parallel.

The main result of this paper proves \cite[Conjecture~2]{AMOS} for surfaces with flat
normal bundle.

\begin{theorem}\label{thm:main}
Let $\varphi\colon M^2\to\Sph^4$ be a proper biharmonic surface with flat normal
bundle. Then, the image $\varphi(M)$ lies minimally in a small hypersphere $\Sph^3(2)$
of $\Sph^4$.
\end{theorem}

Theorem~\ref{thm:main} is local and requires no completeness or compactness
assumption. Conversely, a minimal surface of a small hypersphere $\Sph^3(2)$ is a
proper biharmonic surface of $\Sph^4$ and has flat normal bundle, so
Theorem~\ref{thm:main} characterizes these surfaces.

Since any rotational surface in $\Sph^4$ has flat normal bundle, we obtain the full
classification of biharmonic rotational surfaces in $\Sph^4$.

\begin{corollary}\label{cor:rotational}
A rotational surface $\varphi\colon M^2\to\Sph^4$ is proper biharmonic if and only if
$\varphi(M)$ lies minimally in a small hypersphere $\Sph^3(2)$ of $\Sph^4$.
\end{corollary}

We note that there are surfaces in $\Sph^4$ with flat normal bundle which lie minimally in
$\Sph^3(2)$ and are not rotational; for instance, Lawson's minimal surface $\tau_{2,1}$
of $\Sph^3(2)$ (see \cite{Lawson}). We also note that Corollary \ref{cor:rotational} can
be seen as the analogue of Maeta's result on biharmonic rotational surfaces in $\mathbb R^4$
(see \cite{Maeta2026}).

By \cite{BO}, Theorem~\ref{thm:main} amounts to showing that the mean curvature function
$f$ of a proper biharmonic surface with flat normal bundle in $\Sph^4$ is constant.
Some cases of this statement are already known. Pseudo-umbilical biconservative surfaces
in four-dimensional space forms are CMC (see \cite{CMO,Dimitric}, and
Lemma~\ref{lem:pseudo} below). The PNMC case is also solved: the PNMC biconservative
surfaces in $\Sph^4$ were classified by Nistor, Oniciuc, Turgay and Ye\u{g}in \c{S}en
\cite{NOTY}, and it follows from their classification that there are no PNMC proper
biharmonic surfaces in $\Sph^4$ with non-constant mean curvature (see also \cite{BalmusMontaldoOniciuc2013}); the analogous results
in $\R^4$ and in the hyperbolic space $\mathbb H^4$ were obtained in \cite{Chen2019, YT} and
\cite{NR}, respectively (see Proposition~\ref{prop:PNMC}). Therefore, the core of Theorem~\ref{thm:main} is the
following non-existence result, which is proved in Section~\ref{sec:conj}.

\begin{theorem}\label{thm:nonexistence}
There are no biharmonic surfaces with flat normal bundle in $\Sph^4$ such that, at any
point, $H\neq0$, $\grad f\ne0$, $A_{E_3}\neq f\Id$ and $\nabla^\perp E_3\neq0$.
\end{theorem}

The surfaces of Theorem~\ref{thm:nonexistence} were studied in \cite{AMOS}, in any
four-dimensional space form $N^4(\eps)$, under two additional hypotheses: the surface is
a \emph{Weingarten surface} (W-surface), that is, $f$ and the Gaussian curvature $K$
are functionally dependent, and $3f^2+K-\eps\neq0$ at any point. Under these hypotheses
one has $\langle\nabla^\perp_{E_1}E_3,E_4\rangle\langle\nabla^\perp_{E_2}E_3,E_4\rangle=0$
for a suitable adapted frame, the functions which describe the surface depend on only
one variable, and the biconservative surfaces are characterized by a first order system
of ordinary differential equations. In the biharmonic case, this system is completed by
the normal part of the biharmonic equation, and becomes a first order system of ODEs
with a constraint \cite[Theorem~4.1]{AMOS}. Inspired by the analysis of several explicit
families, we conjectured in \cite[Conjecture~1]{AMOS} that the non-CMC, non-PNMC
biharmonic W-surfaces with flat normal bundle and $3f^2+K-\eps\neq0$ do not exist, and we showed that this would imply \cite[Conjecture~2]{AMOS} for this type of surfaces. 
Theorem~\ref{thm:nonexistence} proves \cite[Conjecture~1]{AMOS} for
$\eps=1$ and, at the same time, removes the Weingarten hypothesis and the hypothesis
$3f^2+K-\eps\neq0$.

The proof of Theorem~\ref{thm:nonexistence} consists of three steps. The first two of
them hold in any space form $N^4(\eps)$, and only the last one uses $\eps=1$.

\emph{The gradient of $f$ is a principal direction.} In \cite{AMOS}, the Weingarten hypothesis and
the condition $3f^2+K-\eps\neq0$ were used precisely to show that $\grad f$ is a principal
direction of the shape operators. For biharmonic surfaces this property holds without
these hypotheses (Proposition~\ref{prop:principal}). For $\eps\leq0$ this is a result of
Fu, Hong and Yang \cite[Lemma~4.3]{FHY}, and we observe that their argument holds for
any $\eps$. The idea is to normalize all the quantities by $f$. Since the normal bundle is
flat, the connection form of the normal bundle is closed, and, together with the
biconservativity, this implies that the function
\[
 \beta=\frac{\left(k_1-k_2\right)\left(3f^2+K-\eps\right)}{2\alpha\left(f^2-K+\eps\right)}
\]
depends, locally, only on $f$ (Lemma~\ref{lem:beta}); here $k_1,k_2$ and $\alpha,-\alpha$
are the eigenvalues of the shape operators $A_{E_3}$ and $A_{E_4}$, where $E_4$ is a unit
normal vector field orthogonal to $E_3$. If $\grad f$ were not a
principal direction on an open set, then $f$ and $z=2\alpha/(k_1-k_2)$ would be local
coordinates, and the compatibility conditions for the second derivatives of $f$, together
with the normal part of the biharmonic equation, would give a linear system for the
functions $|\grad f|^2/f^4$ and $(k_1-k_2)\bigl((E_1f)^2-(E_2f)^2\bigr)/(2f^5)$ whose
coefficients are polynomials in $z$ depending only on $f$. Hence, along the level sets of
$f$, these two functions are rational functions of $z$. A differential relation along the
level sets of $f$, evaluated at $z=\infty$, then forces $f^2=\eps$, which is a
contradiction.

\emph{Reduction to ODEs.} Once $\grad f$ is a principal direction, a short argument which
uses only the biconservativity, the Codazzi and the Ricci equations shows that all the
functions which describe the surface depend on one variable $s$, the arc length along
the integral curves of $\grad f/|\grad f|$ (Proposition~\ref{prop:reduction}). In
particular, the surface is, automatically, a W-surface, and the biharmonic equation is
equivalent to system (4.1) of \cite{AMOS} (Theorem~\ref{thm:41}). Thus, the normal part of
the biharmonic equation compensates for the two additional hypotheses of \cite{AMOS}.

\emph{Elimination.} For $\eps=1$, we normalize the unknown functions of the ODE system
by the mean curvature function and change the parameter $s$ into $\tau$, with
$d\tau/ds=f$. The system becomes a polynomial vector field in six variables, one of which
is $e=1/f^2$, together with a polynomial constraint $C=0$. The constraint and its first
derivative along the flow determine two of the variables, and the next four derivatives
give four polynomial equations in the remaining four variables. An
elimination argument based on resultants shows that, on the common zero set of these
equations, $e$ takes only finitely many values. Hence $f$ is constant, which is a
contradiction (Proposition~\ref{prop:finite}).

Some of the identities in Sections~\ref{sec:principal} and \ref{sec:conj} involve long
expressions, and the resultants in Section~\ref{sec:conj} are polynomials with several
thousand terms. All of them were obtained and checked with exact symbolic computations; no numerical approximation is used.

The paper is organized as follows. Section~\ref{sec:prelim} recalls the biharmonic
equation, the adapted frame of a surface with flat normal bundle and the corresponding
fundamental equations. Section~\ref{sec:principal} shows that $\grad f$ is a principal
direction, and Section~\ref{sec:ODE} reduces the problem to a system of ODEs.
Theorem~\ref{thm:nonexistence} is proved in Section~\ref{sec:conj}. Finally, in
Section~\ref{sec:proof}, we recall the PNMC case and prove Theorem~\ref{thm:main}. In
Appendix~\ref{app:rotational}, we describe the geometry of the surfaces which satisfy the
hypotheses of Proposition~\ref{prop:41minimal}, where we slightly relax the hypothesis of Theorem \ref{thm:41} without changing the conclusion. In particular, for $\eps=1$, we show that such surfaces are rotational and we characterize the
biharmonic ones by means of their profile curves. In Appendix \ref{app:code}, we give the Mathematica code which reproduces the
	algebraic elimination in the proof of Proposition~\ref{prop:finite}.
	
\medskip
\noindent \emph{Acknowledgments.} Thanks are due to Professor Yu Fu for valuable comments and suggestions.

\medskip
\noindent\emph{Use of AI tools.} This paper was prepared in interaction with the AI
assistant Claude Opus~5.5 (Anthropic), which was used to check computations, to suggest
expository improvements and to assist in editing the manuscript. All the mathematical
content, including the statements and proofs, has been deduced and verified by the
authors, who take full responsibility for it.

\section{Preliminaries}\label{sec:prelim}

All manifolds, maps and tensor fields are smooth, all manifolds are connected and all
immersions are isometric. The arguments are local, and, when no ambiguity can arise, we
identify tangent vectors with their images under the differential of the immersion.
We denote by $N^4(\eps)$ the four-dimensional space form of constant sectional curvature
$\eps$. Although our main interest is the case $\eps=1$, the results of
Sections~\ref{sec:prelim}--\ref{sec:ODE} hold for any $\eps\in\R$, and we keep $\eps$ to
make this apparent.

\subsection{Biharmonic surfaces in space forms}
Our curvature convention is $R(X,Y)Z=\nabla_X\nabla_YZ-\nabla_Y\nabla_XZ-\nabla_{[X,Y]}Z$,
and the same convention is used for the curvature $R^\perp$ of the normal bundle. For a
metric connection $\nabla^{\mathcal E}$ on a vector bundle, the rough Laplacian is
\[
 \Delta^{\mathcal E}V=-\sum_{i=1}^m
 \left(\nabla^{\mathcal E}_{e_i}\nabla^{\mathcal E}_{e_i}V
 -\nabla^{\mathcal E}_{\nabla_{e_i}e_i}V\right),
\]
where $\{e_i\}$ is a local orthonormal frame of $TM$.

Let $\varphi\colon M^2\to N^4(\eps)$ be a surface. The second fundamental form $B$, the
shape operator $A_\xi$ associated with a normal vector field $\xi$, and the normal
connection $\nabla^\perp$ are defined by the Gauss and Weingarten formulas
\[
 \nabla^\varphi_XY=\nabla_XY+B(X,Y),\qquad
 \nabla^\varphi_X\xi=-A_\xi X+\nabla^\perp_X\xi,\qquad
 \langle A_\xi X,Y\rangle=\langle B(X,Y),\xi\rangle,
\]
and $H=\frac12\tr B$. We recall the Gauss, Codazzi and Ricci equations: for tangent
vector fields $X,Y,Z,W$ and normal vector fields $\xi,\eta$,
\begin{gather}
 \begin{aligned}
 \langle R(X,Y)Z,W\rangle={}&\eps\left(\langle Y,Z\rangle\langle X,W\rangle
 -\langle X,Z\rangle\langle Y,W\rangle\right)\\
 &+\langle B(Y,Z),B(X,W)\rangle-\langle B(X,Z),B(Y,W)\rangle,
 \end{aligned}\label{eq:generalgauss}\\
 (\nabla_XA_\xi)Y-(\nabla_YA_\xi)X=A_{\nabla^\perp_X\xi}Y-A_{\nabla^\perp_Y\xi}X,
 \label{eq:generalcodazzi}\\
 \langle R^\perp(X,Y)\xi,\eta\rangle=\langle[A_\xi,A_\eta]X,Y\rangle .
 \label{eq:generalricci}
\end{gather}
The normal bundle is \emph{flat} when $R^\perp=0$; by \eqref{eq:generalricci}, this is
equivalent to the commutativity of all the shape operators.

The decomposition of the bitension field into its normal and tangent components gives
the following characterization (see \cite{Chen,Oniciuc}).

\begin{proposition}\label{prop:biharmonic}
A surface $\varphi\colon M^2\to N^4(\eps)$ is biharmonic if and only if
\begin{equation}\label{eq:normal}
 \Delta^\perp H+\tr B(\cdot,A_H(\cdot))-2\eps H=0
\end{equation}
and
\begin{equation}\label{eq:tangent}
 2\tr A_{\nabla^\perp_{(\cdot)}H}(\cdot)+\grad|H|^2=0 .
\end{equation}
\end{proposition}

The surface $\varphi$ is called \emph{biconservative} if \eqref{eq:tangent} holds. Thus,
every biharmonic surface is biconservative.

On the open set where $H\neq0$, we write $H=fE_3$, where $f=|H|>0$ is the mean curvature
function and $E_3$ is a unit normal vector field. Following \cite{AMOS}, a surface is
called \emph{non-CMC} if $H\neq0$, $\grad f\neq0$ and $A_3\neq f\Id$ at any point, where
$A_3=A_{E_3}$. The last condition means that the surface has no pseudo-umbilical points.
It is natural in view of the following lemma.

\begin{lemma}[\cite{CMO,Dimitric}]\label{lem:pseudo}
Let $\varphi\colon M^2\to N^4(\eps)$ be a biconservative surface such that $H\neq0$ and
$A_3=f\Id$ on an open subset $U$. Then $f$ is constant on $U$.
\end{lemma}

\begin{proof}
Let $E_4$ be a unit normal vector field orthogonal to $E_3$, $A_4=A_{E_4}$, and let
$\omega(X)=\langle\nabla^\perp_XE_3,E_4\rangle$, so that $\nabla^\perp_XE_3=\omega(X)E_4$.
Denote by $\omega^\sharp$ the vector field dual to $\omega$. Since $\tr A_4=2\langle
H,E_4\rangle=0$, $A_4$ is a traceless symmetric operator. As $\nabla A_3=df\otimes\Id$,
the Codazzi equation \eqref{eq:generalcodazzi} for $\xi=E_3$ gives
\[
 X(f)Y-Y(f)X=\omega(X)A_4Y-\omega(Y)A_4X .
\]
Taking $X=e_1$, $Y=e_2$ for a local orthonormal frame $\{e_1,e_2\}$, and using
$\langle A_4e_1,e_1\rangle=-\langle A_4e_2,e_2\rangle$, we obtain
$\grad f=-A_4\omega^\sharp$. On the other hand,
$\nabla^\perp_{e_i}H=e_i(f)E_3+f\omega(e_i)E_4$, so that
$\tr A_{\nabla^\perp_{(\cdot)}H}(\cdot)=f\grad f+fA_4\omega^\sharp$, and
\eqref{eq:tangent} becomes $2\grad f+A_4\omega^\sharp=0$. Therefore $\grad f=0$.
\end{proof}

\subsection{The adapted frame}
Let $\varphi\colon M^2\to N^4(\eps)$ be a non-CMC surface with flat normal bundle. Since
the normal bundle has rank two, there exists, locally, a unit normal vector field $E_4$
orthogonal to $E_3$, unique up to sign. We denote $A_4=A_{E_4}$. The following
proposition collects the fundamental equations in a suitable frame; it is
\cite[Propositions~3.1 and 3.2]{AMOS} and the computation in the proof of
\cite[Theorem~4.1]{AMOS}, where the W-hypothesis is not used.

\begin{proposition}\label{prop:frame}
Let $\varphi\colon M^2\to N^4(\eps)$ be a non-CMC surface with flat normal bundle. Then,
locally, there exists an orthonormal frame field $\{E_1,E_2\}$ of $TM$ such that
\begin{equation}\label{eq:frame}
 A_3E_1=k_1E_1,\quad A_3E_2=k_2E_2,\quad A_4E_1=\alpha E_1,\quad A_4E_2=-\alpha E_2,
\end{equation}
where $k_1+k_2=2f$ and $k_1\neq k_2$ at any point. Write
\begin{gather}
 \nabla_{E_1}E_1=-a_1E_2,\quad \nabla_{E_1}E_2=a_1E_1,\quad \nabla_{E_2}E_1=-a_2E_2,\quad
 \nabla_{E_2}E_2=a_2E_1,\label{eq:connection}\\
 \nabla^\perp_{E_i}E_3=b_iE_4,\qquad \nabla^\perp_{E_i}E_4=-b_iE_3,\qquad i=1,2 .
 \label{eq:normalconnection}
\end{gather}
Then $[E_1,E_2]=a_1E_1+a_2E_2$, and the following hold: the Ricci equation
\begin{equation}\label{eq:ricci}
 E_1(b_2)-E_2(b_1)=a_1b_1+a_2b_2,
\end{equation}
the Gauss equation
\begin{equation}\label{eq:gauss}
 K=E_1(a_2)-E_2(a_1)-a_1^2-a_2^2=\eps+k_1k_2-\alpha^2,
\end{equation}
and the Codazzi equations
\begin{equation}\label{eq:codazziAMOS}
\begin{gathered}
 E_1(k_2)=a_2(k_2-k_1)-\alpha b_1,\qquad E_2(k_1)=a_1(k_2-k_1)+\alpha b_2,\\
 E_1(\alpha)=2\alpha a_2+k_2b_1,\qquad E_2(\alpha)=-2\alpha a_1-k_1b_2 .
\end{gathered}
\end{equation}
Moreover, the surface is biconservative if and only if
\begin{equation}\label{eq:bicons}
 (k_1+f)E_1(f)=-f\alpha b_1,\qquad (k_2+f)E_2(f)=f\alpha b_2 ,
\end{equation}
and it is biharmonic if and only if \eqref{eq:bicons} holds and
\begin{equation}\label{eq:normalframe}
\left\{
\begin{aligned}
&-E_1(E_1(f))-E_2(E_2(f))+f\left(b_1^2+b_2^2+k_1^2+k_2^2-2\eps\right)-a_1E_2(f)+a_2E_1(f)=0\\
&2\left(b_1E_1(f)+b_2E_2(f)\right)+f\left(E_1(b_1)+E_2(b_2)\right)+f\left(a_1b_2-a_2b_1\right)+f\alpha\left(k_2-k_1\right)=0
\end{aligned}\right. .
\end{equation}
\end{proposition}

\begin{proof}
Since the normal bundle is flat, \eqref{eq:generalricci} gives $[A_3,A_4]=0$. As
$A_3\neq f\Id$ and $\tr A_3=2f$, the operator $A_3$ has two distinct eigenvalues $k_1$,
$k_2$ at any point, so they are smooth and there is a local orthonormal frame
$\{E_1,E_2\}$ of eigenvectors. Since $A_4$ commutes with $A_3$, it is diagonal in this
frame, and since $\tr A_4=2\langle H,E_4\rangle=0$, we obtain \eqref{eq:frame}. The
expression of $[E_1,E_2]$ follows from \eqref{eq:connection}, and \eqref{eq:ricci} is
$R^\perp(E_1,E_2)E_3=0$. The first expression of $K$ in \eqref{eq:gauss} is
$\langle R(E_1,E_2)E_2,E_1\rangle$ computed with \eqref{eq:connection}, and the second
one follows from \eqref{eq:generalgauss}, since $B(E_1,E_1)=k_1E_3+\alpha E_4$,
$B(E_1,E_2)=0$ and $B(E_2,E_2)=k_2E_3-\alpha E_4$. The equations
\eqref{eq:codazziAMOS} are the $E_1$ and $E_2$ components of \eqref{eq:generalcodazzi}
with $X=E_1$, $Y=E_2$ and $\xi=E_3$, $\xi=E_4$.

Since $\nabla^\perp_{E_i}H=E_i(f)E_3+fb_iE_4$, we have
\[
 \tr A_{\nabla^\perp_{(\cdot)}H}(\cdot)=\left(k_1E_1(f)+f\alpha b_1\right)E_1
 +\left(k_2E_2(f)-f\alpha b_2\right)E_2,
\]
and \eqref{eq:tangent} is equivalent to \eqref{eq:bicons}. Finally, a direct computation
gives
\begin{align*}
 \Delta^\perp H={}&\left(-E_1(E_1(f))-E_2(E_2(f))+f\left(b_1^2+b_2^2\right)-a_1E_2(f)+a_2E_1(f)\right)E_3\\
 &-\left(2\left(b_1E_1(f)+b_2E_2(f)\right)+f\left(E_1(b_1)+E_2(b_2)\right)+f\left(a_1b_2-a_2b_1\right)\right)E_4,\\
 \tr B(\cdot,A_H(\cdot))={}&f\left(k_1^2+k_2^2\right)E_3+f\alpha\left(k_1-k_2\right)E_4,
\end{align*}
so \eqref{eq:normal} is equivalent to \eqref{eq:normalframe}.
\end{proof}

\begin{remark}\label{rem:symmetry}
Interchanging $E_1$ and $E_2$ and replacing $E_4$ by $-E_4$, the functions
\[
 (k_1,k_2,\alpha,a_1,a_2,b_1,b_2)\quad\text{become}\quad(k_2,k_1,\alpha,-a_2,-a_1,-b_2,-b_1),
\]
and the set of equations \eqref{eq:ricci}--\eqref{eq:normalframe} is unchanged.
\end{remark}

We first exclude the case $A_4=0$ for non-PNMC surfaces.

\begin{lemma}\label{lem:A4}
Let $\varphi\colon M^2\to N^4(\eps)$ be a non-CMC biconservative surface with flat
normal bundle. If $A_4=0$ on an open subset $U$ of $M$, then there exists a non-empty
open subset of $U$ on which $\nabla^\perp E_3=0$.
\end{lemma}

\begin{proof}
On $U$ we have $\alpha=0$, and \eqref{eq:bicons} gives $(k_1+f)E_1(f)=(k_2+f)E_2(f)=0$.
Since $\grad f\neq0$, eventually restricting $U$ and using Remark~\ref{rem:symmetry}, we
can assume that $E_1(f)\neq0$ at any point of $U$. Then $k_1=-f$ and $k_2=3f$, and the
Codazzi equations \eqref{eq:codazziAMOS} for $\alpha$ give $k_2b_1=0$ and $k_1b_2=0$.
Thus $b_1=b_2=0$, that is, $\nabla^\perp E_3=0$ on $U$.
\end{proof}

\section{The gradient of the mean curvature is a principal direction}\label{sec:principal}

In this section, $\varphi\colon M^2\to N^4(\eps)$ is a non-CMC biharmonic surface with
flat normal bundle, $\eps\in\R$ is arbitrary, and we work on an open subset on which
the frame of Proposition~\ref{prop:frame} is defined and $\alpha\neq0$ at any point. We
show that $\grad f$ is a principal direction of $A_3$ and $A_4$, following the proof of
\cite[Lemma~4.3]{FHY}. We use the auxiliary functions of \cite{FHY}, written in our
notation:
\begin{equation}\label{eq:aux}
 \lambda=\frac{k_1-k_2}{2f},\qquad \mu=\frac{\alpha}{f},\qquad u=\frac{E_1(f)}{f},\qquad
 v=\frac{E_2(f)}{f},\qquad r=\lambda^2+\mu^2 .
\end{equation}
Then $k_1=f(1+\lambda)$, $k_2=f(1-\lambda)$, $\alpha=f\mu$, $\lambda\mu\neq0$,
$\grad f=f(uE_1+vE_2)$, and, by \eqref{eq:gauss},
\begin{equation}\label{eq:DeltaL}
 3f^2+K-\eps=f^2(4-r),\qquad f^2-K+\eps=f^2r .
\end{equation}
The relations \eqref{eq:bicons} become
\begin{equation}\label{eq:b}
 b_1=-\frac{2+\lambda}{\mu}\,u,\qquad b_2=\frac{2-\lambda}{\mu}\,v,
\end{equation}
and, replacing \eqref{eq:b} in the Codazzi equations \eqref{eq:codazziAMOS}, we obtain
\begin{equation}\label{eq:codazzi}
\left\{
\begin{aligned}
 E_1(\lambda)&=-(1+2\lambda)u+2\lambda a_2, &
 E_2(\lambda)&=(1-2\lambda)v-2\lambda a_1,\\
 E_1(\mu)&=-\frac{2-\lambda-\lambda^2+\mu^2}{\mu}\,u+2\mu a_2, &
 E_2(\mu)&=-\frac{2+\lambda-\lambda^2+\mu^2}{\mu}\,v-2\mu a_1 .
\end{aligned}\right.
\end{equation}
The vector field $\grad f$ is a principal direction of $A_3$ and $A_4$ at a point if and
only if $uv=0$ at that point.

Since the normal bundle is flat, by \eqref{eq:ricci} the $1$-form
$\omega=b_1\omega^1+b_2\omega^2$, where $\{\omega^1,\omega^2\}$ is the dual coframe of
$\{E_1,E_2\}$, is closed. Thus, locally, $\omega=d\theta$ for some function $\theta$,
that is, $b_i=E_i(\theta)$.

Assume now that $uv\neq0$ at any point of an open subset. Following \cite{FHY}, we put
\begin{equation}\label{eq:PQ}
 P=-\frac{a_1}{v},\qquad Q=-\frac{a_2}{u},\qquad z=\frac{\mu}{\lambda},\qquad
 \tE_1=\frac{1}{u}E_1,\qquad \tE_2=\frac{1}{v}E_2 .
\end{equation}
Then $\tE_1(f)=\tE_2(f)=f$.

\begin{lemma}[{\cite[Lemma~4.1]{FHY}}]\label{lem:beta}
On an open subset where $uv\neq0$, the function
\begin{equation}\label{eq:beta}
 \beta=\frac{\lambda(4-r)}{\mu r}
 =\frac{\left(k_1-k_2\right)\left(3f^2+K-\eps\right)}{2\alpha\left(f^2-K+\eps\right)}
\end{equation}
depends, locally, only on $f$.
\end{lemma}

\begin{proof}
From \eqref{eq:b} we have $\tE_1(\theta)=-(2+\lambda)/\mu$ and
$\tE_2(\theta)=(2-\lambda)/\mu$. From \eqref{eq:codazzi}, the terms containing $a_1$ and
$a_2$ cancel in the derivatives of $z$, and we get
\[
 \tE_1(z)=\frac{(\lambda+1)r-2\lambda}{\mu\lambda^2},\qquad
 \tE_2(z)=\frac{(\lambda-1)r-2\lambda}{\mu\lambda^2}.
\]
Since $1+z^2=r/\lambda^2$, it follows that
\[
\begin{aligned}
 \tE_1\left(\theta+2\arctan z\right)&=\frac{-(2+\lambda)r+2(\lambda+1)r-4\lambda}{\mu r}=-\beta,\\
 \tE_2\left(\theta+2\arctan z\right)&=\frac{(2-\lambda)r+2(\lambda-1)r-4\lambda}{\mu r}=-\beta .
\end{aligned}
\]
As $\tE_1(\log f)=\tE_2(\log f)=1$, we get $d\left(\theta+2\arctan z\right)=-\beta\,
d\log f$. Taking the exterior derivative, we obtain $d\beta\wedge df=0$ and, since
$df\neq0$, $\beta$ is, locally, a function of $f$.
\end{proof}

We denote
\begin{equation}\label{eq:gamma}
 \gamma=f\frac{d\beta}{df},
\end{equation}
and by $\gamma'$ the derivative of $\gamma$ with respect to $f$. Then
$\tE_1(\beta)=\tE_2(\beta)=\gamma$ and $\tE_1(\gamma)=\tE_2(\gamma)=f\gamma'$.

\begin{lemma}[{\cite[(4.19) and Lemma~4.2]{FHY}}]\label{lem:derivatives}
On an open subset where $uv\neq0$, the following hold:
\begin{enumerate}[label=\alph*)]
\item the functions $P$ and $Q$ are given by
\begin{align}
 P+Q&=\frac{r(4-r)}{8\lambda\mu^2},\label{eq:P+Q}\\
 P-Q&=\frac{r+2}{r}-\frac{(4-r)(r-2)}{8\mu^2}-\frac{r\mu}{8\lambda}\,\gamma ;\label{eq:P-Q}
\end{align}
\item the derivatives of $u$ and $v$ are given by
\begin{align}
 E_2(u)&=\left(2P-Q-\frac{3\mu^2+4-\lambda^2}{2\mu^2}\right)uv,\qquad
 E_1(v)=\left(P-2Q-\frac{3\mu^2+4-\lambda^2}{2\mu^2}\right)uv,\label{eq:mixed}\\
 E_1(u)&=\left(\frac{(\lambda+2)^2}{4\mu^2}-\frac74-2Q\right)u^2
 +\left(\frac{3(\lambda-2)^2}{4\mu^2}+\frac34-P\right)v^2\notag\\
 &\quad+\left(1+\lambda^2-\frac{\lambda(1+r)}{2}\right)f^2-\frac{2-\lambda}{2}\,\eps,\label{eq:E1u}\\
 E_2(v)&=\left(\frac{3(\lambda+2)^2}{4\mu^2}+\frac34+Q\right)u^2
 +\left(\frac{(\lambda-2)^2}{4\mu^2}-\frac74+2P\right)v^2\notag\\
 &\quad+\left(1+\lambda^2+\frac{\lambda(1+r)}{2}\right)f^2-\frac{2+\lambda}{2}\,\eps;\label{eq:E2v}
\end{align}
\item the derivative $\gamma'$ is determined by
\begin{multline}\label{eq:gammaprime}
 \frac{\mu r\left(u^2+v^2\right)}{16\lambda}\,f\gamma'
 =\left(P_\lambda \frac{E_2(\lambda)}{v}+P_\mu \frac{E_2(\mu)}{v}-P^2\right)v^2
 -\left(Q_\lambda \frac{E_1(\lambda)}{u}+Q_\mu \frac{E_1(\mu)}{u}+Q^2\right)u^2\\
 +PE_2(v)-QE_1(u)-\eps-f^2(1-r),
\end{multline}
where $P_\lambda$, $P_\mu$, $Q_\lambda$ and $Q_\mu$ denote the partial derivatives of
$P$ and $Q$, given by \eqref{eq:P+Q}--\eqref{eq:P-Q}, as functions of $(\lambda,\mu)$,
with $\gamma$ fixed.
\end{enumerate}
\end{lemma}

\begin{proof}
a) By \eqref{eq:codazzi}, $E_1(\beta)=\beta_\lambda E_1(\lambda)+\beta_\mu E_1(\mu)$ is a
linear function of $u$ and $a_2$, in which the coefficient of $a_2$ is
$2(\lambda\beta_\lambda+\mu\beta_\mu)=-16\lambda/(\mu r)\neq0$. Since $\beta$ depends
only on $f$, we have $E_1(\beta)=\gamma u$. Replacing $a_2=-Qu$ and dividing by $u$, we
obtain $Q$ as a function of $\lambda$, $\mu$ and $\gamma$. In the same way, the relation
$E_2(\beta)=\gamma v$ determines $P$, and a direct computation gives
\eqref{eq:P+Q}--\eqref{eq:P-Q}.

b) The commutator identity $E_1(E_2(f))-E_2(E_1(f))=a_1E_1(f)+a_2E_2(f)$ gives
\[
 E_1(v)-E_2(u)=-(P+Q)uv .
\]
Replacing \eqref{eq:b} in the Ricci equation \eqref{eq:ricci} and using
\eqref{eq:codazzi}, we obtain
\[
 (2-\lambda)E_1(v)+(2+\lambda)E_2(u)=\left(\lambda(P+Q)+6(P-Q)-6+\frac{2(\lambda^2-4)}{\mu^2}\right)uv .
\]
Solving this linear system, we get \eqref{eq:mixed}. Similarly, replacing \eqref{eq:b},
$a_1=-Pv$ and $a_2=-Qu$ in \eqref{eq:normalframe}, the two components of the normal
part of the biharmonic equation become
\[
\begin{aligned}
 E_1(u)+E_2(v)={}&\left(\frac{(\lambda+2)^2}{\mu^2}-1-Q\right)u^2
 +\left(\frac{(\lambda-2)^2}{\mu^2}-1+P\right)v^2\\
 &+2(1+\lambda^2)f^2-2\eps,\\
 (\lambda+2)E_1(u)+(\lambda-2)E_2(v)={}&-\left(Q(\lambda+6)+\lambda+5-\frac{(\lambda-1)(\lambda+2)^2}{\mu^2}\right)u^2\\
 &+\left(P(\lambda-6)+5-\lambda+\frac{(\lambda+1)(\lambda-2)^2}{\mu^2}\right)v^2\\
 &-2\lambda\mu^2f^2 .
\end{aligned}
\]
The determinant of this linear system in $E_1(u)$ and $E_2(v)$ is $-4$, and its solution is
\eqref{eq:E1u}--\eqref{eq:E2v}.

c) We replace $a_1=-Pv$ and $a_2=-Qu$ in the Gauss equation
$E_1(a_2)-E_2(a_1)-a_1^2-a_2^2=\eps+f^2(1-r)$. By \eqref{eq:P-Q}, the derivatives of $P$
and $Q$ with respect to $\gamma$ are $-r\mu/(16\lambda)$ and $r\mu/(16\lambda)$,
respectively, and $E_1(\gamma)=f\gamma'u$, $E_2(\gamma)=f\gamma'v$. Thus the terms
containing $\gamma'$ are $-r\mu f\gamma'(u^2+v^2)/(16\lambda)$, and we obtain
\eqref{eq:gammaprime}.
\end{proof}

We note that the sign of $\eps$ plays no role in the computations above.

\begin{remark}\label{rem:PQAMOS}
Substituting \eqref{eq:b}, $a_1=-Pv$ and $a_2=-Qu$, and using \eqref{eq:DeltaL}, one
checks that, at the points where $(k_1+f)(k_2+f)\neq0$, relation \eqref{eq:P+Q} is
equivalent to
\[
 4f\alpha\left(k_1-k_2\right)\left(a_1b_1\left(k_2+f\right)-a_2b_2\left(k_1+f\right)\right)
 +b_1b_2\left(3f^2+K-\eps\right)\left(f^2-K+\eps\right)=0 .
\]
For W-surfaces, one also has $a_1E_1(f)+a_2E_2(f)=0$ \cite[(3.23)]{AMOS}, and then this
relation gives $b_1b_2\left(3f^2+K-\eps\right)=0$, which is the conclusion of
\cite[Lemma~3.3]{AMOS}.
\end{remark}

\begin{proposition}\label{prop:principal}
Let $\varphi\colon M^2\to N^4(\eps)$ be a non-CMC biharmonic surface with flat normal
bundle. Assume that $A_4\neq0$ at any point. Then, $\grad f$ is a principal direction of
$A_3$ and $A_4$ at any point.
\end{proposition}

\begin{proof}
Assume, by way of contradiction, that $uv\neq0$ at any point of an open subset. Let
\[
 \Gf=\frac{u^2+v^2}{f^2}=\frac{|\grad f|^2}{f^4},\qquad
 \Zf=\frac{\lambda\left(u^2-v^2\right)}{f^2}.
\]
Since $\tE_1(f)=\tE_2(f)=f$ and $\left(\tE_1-\tE_2\right)(z)=2r/\left(\mu\lambda^2\right)\neq0$,
the functions $(f,z)$ are local coordinates. By Lemma~\ref{lem:beta}, the functions
$\beta$, $\gamma$ and $\eps/f^2$ are constant when $f$ is fixed. Moreover,
$\beta z=(4-r)/r$, so that
\[
 r=\frac{4}{1+\beta z},\qquad \lambda^2=\frac{r}{1+z^2},\qquad \mu=\lambda z .
\]

We consider the compatibility conditions for $u$ and $v$, that is, the commutator
identities
\[
\begin{gathered}
 E_2(E_1(u))-E_1(E_2(u))-PvE_1(u)-QuE_2(u)=0,\\
 E_2(E_1(v))-E_1(E_2(v))-PvE_1(v)-QuE_2(v)=0,
\end{gathered}
\]
which we compute using Lemma~\ref{lem:derivatives} and \eqref{eq:codazzi}. In these
relations $\gamma'$ only appears through the term $f\gamma'(u^2+v^2)$, which we replace
using \eqref{eq:gammaprime}, and the first relation is divisible by $v$, the second one
by $u$. Dividing the first relation by $vf^2$ and the second one by $uf^2$, adding them
and multiplying the result by $-2\lambda z^3(1+\beta z)$, we obtain
{\small
\begin{multline}\label{eq:GZ1}
\left[\beta^2z\left(3z^4+10z^2+11\right)+2\beta\left(z^4+7\right)+4\gamma z^2-12z\right]\Gf\\
-\left(1+z^2\right)\left[\beta^2z\left(z^2-9\right)+\beta\left(9z^2-7\right)+2\gamma z^2+6z\right]\Zf
=4\beta z^2\left(z^2+5\right)\frac{\eps}{f^2}\\
-\frac{4z^2}{\left(1+z^2\right)\left(1+\beta z\right)^2}\Bigl[\beta^3z^2\left(1+z^2\right)\left(z^2+5\right)\\
+2\beta^2z\left(z^4+2z^2+15\right)+\beta\left(z^4-22z^2+17\right)-8\gamma z^2-12z\Bigr].
\end{multline}}
Subtracting them and multiplying the result by $4z^3(1+\beta z)/(1+z^2)$, we obtain
{\small
\begin{multline}\label{eq:GZ2}
\left[3\beta^3z^2\left(z^2+3\right)+\beta^2z\left(13z^2+21\right)+6\beta\gamma z^3+2\beta\left(6z^2+7\right)+8\gamma z^2\right]\Gf
+7\beta\left(1+z^2\right)\left(1+\beta z\right)^2\Zf\\
=\frac{4z^2}{1+z^2}\left[\beta^2z\left(z^2+3\right)+\beta\left(8z^2+5\right)+2\gamma z^2+5z\right]\frac{\eps}{f^2}\\
-\frac{4z^2}{\left(1+z^2\right)\left(1+\beta z\right)}\Bigl[\beta^3z^2\left(z^2+3\right)
+3\beta^2z\left(3z^2+8\right)+\beta\left(2\gamma z^3+13z^2+17\right)+2\gamma z^2+9z\Bigr].
\end{multline}}
All the quantities in \eqref{eq:GZ1}--\eqref{eq:GZ2}, except $\Gf$ and $\Zf$, are
functions of $z$, $\beta$, $\gamma$ and $\eps/f^2$, that is, for fixed $f$, they are
rational functions of $z$.

If $\beta=0$ on an open subset, then $\gamma=0$, and \eqref{eq:GZ2} reduces to
$z^3\left(5\eps-9f^2\right)=0$. Since $z\neq0$, $f$ is constant, which is a
contradiction. Thus, eventually restricting the open subset, we can assume that
$\beta\neq0$ at any point. For fixed $f$, \eqref{eq:GZ1}--\eqref{eq:GZ2} is a linear
system in $\Gf$ and $\Zf$, whose determinant is
\begin{multline*}
 4z^3(1+z^2)\bigl(6\beta^5z^4+3\beta^3\gamma z^3+24\beta^4z^3+13\beta^5z^2+49\beta^3z^2
 +22\beta^2\gamma z^2+3\beta\gamma^2z^2+26\beta^4z\\
 -2\beta^3\gamma z+50\beta^2z+33\beta\gamma z+4\gamma^2z-\beta^5+11\beta^3-4\beta^2\gamma
 +18\beta+12\gamma\bigr),
\end{multline*}
a polynomial in $z$ of degree $9$ with leading coefficient $24\beta^5\neq0$. Thus $\Gf$
and $\Zf$ are rational functions of $z$, and, as $z\to\infty$,
\begin{equation}\label{eq:asympt}
 \Zf=\frac{13\left(\eps/f^2-1\right)}{6\beta^2z^2}+O\left(z^{-3}\right),\qquad
 \Gf=\frac{4\left(\eps/f^2-1\right)}{3\beta z}-\frac{\eps/f^2-1}{6\beta^2z^2}+O\left(z^{-3}\right).
\end{equation}

On the other hand, with $f$ fixed, we have
$\lambda\left(\tE_1-\tE_2\right)=2\left(1+z^2\right)z^{-1}\partial_z$, and, computing
$\lambda\left(\tE_1-\tE_2\right)(\Gf)$ with Lemma~\ref{lem:derivatives}, we obtain
\begin{equation}\label{eq:dGdz}
 \frac{2\left(1+z^2\right)}{z}\frac{\partial\Gf}{\partial z}
 =\left(\frac{r-4}{2}+\frac{r-12}{2z^2}\right)\Gf-\frac{2\left(1+z^2\right)}{z^2}\Zf
 +\frac{2r}{1+z^2}\left(\frac{\eps}{f^2}-r-1\right),
\end{equation}
where $r=4/(1+\beta z)$. Here the terms containing $P$ and $Q$ only appear through
$P+Q$, which is given by \eqref{eq:P+Q}. Both sides of \eqref{eq:dGdz} are rational
functions of $z$ which coincide on an open interval, so they coincide identically.
Substituting \eqref{eq:asympt} in \eqref{eq:dGdz}, the coefficients of $z^{-1}$
coincide, and comparing the coefficients of $z^{-2}$ we get
\[
 \frac{2\left(\eps/f^2-1\right)}{\beta^2}=0,
\]
that is, $f^2=\eps$ on an open subset, which contradicts $\grad f\neq0$.
\end{proof}

\begin{remark}\label{rem:FHY}
The proof of Proposition~\ref{prop:principal} is the proof of \cite[Lemma~4.3]{FHY}. In
\cite{FHY} the condition $\eps\leq0$ is used only at the end, in the form
$\eps/f^2-1<0$; however, $\eps/f^2-1=0$ on an open subset is also impossible, since it
implies that $f$ is constant. In the same way, in the case $\beta=0$, the relation
$5\eps=9f^2$ implies that $f$ is constant. From \eqref{eq:DeltaL}, $\beta=0$
if and only if $3f^2+K-\eps=0$, so the hypothesis $3f^2+K-\eps\neq0$ of \cite{AMOS} is not
needed in Proposition~\ref{prop:principal}.
\end{remark}

\section{Reduction to a system of ODEs}\label{sec:ODE}

\begin{proposition}\label{prop:reduction}
Let $\varphi\colon M^2\to N^4(\eps)$ be a non-CMC biharmonic surface with flat normal
bundle. Assume that $\nabla^\perp E_3\neq0$ at any point. Then, the set
$M'=\{p\in M\mid A_4(p)\neq0\}$ is open and dense in $M$ and, around any point of $M'$,
eventually changing the frame as in Remark~\ref{rem:symmetry}, we have
\begin{equation}\label{eq:340}
 \alpha\neq0,\qquad b_1\neq0,\qquad E_2(f)=0,\qquad b_2=0,\qquad a_1=0
\end{equation}
and
\begin{equation}\label{eq:342}
 E_2(b_1)=E_2(k_1)=E_2(\alpha)=E_2(k_2)=E_2(K)=E_2(a_2)=0 .
\end{equation}
\end{proposition}

\begin{proof}
Clearly, $M'$ is open. By Lemma~\ref{lem:A4}, if the interior of $M\setminus M'$ were
non-empty, there would exist an open subset on which $\nabla^\perp E_3=0$, which
contradicts the hypothesis. Thus, $M'$ is dense.

Let $U\subset M'$ be a connected open subset on which the frame $\{E_1,E_2\}$ is defined.
By Proposition~\ref{prop:principal}, $uv=0$ on $U$. Since $\grad f\neq0$, the sets
$\{u\neq0\}$ and $\{v\neq0\}$ are open, disjoint and cover $U$. Since $U$ is connected,
one of them is empty and, by Remark~\ref{rem:symmetry}, we can assume that $v=0$ on $U$,
that is, $E_2(f)=0$. Then \eqref{eq:b} gives $b_2=0$ and, since $\nabla^\perp E_3\neq0$,
we obtain $b_1\neq0$ at any point of $U$.

Since $E_2(f)=0$, we have $[E_1,E_2](f)=-E_2(fu)=-fE_2(u)$. On the other hand,
$[E_1,E_2]=a_1E_1+a_2E_2$, so $[E_1,E_2](f)=a_1fu$, and we get $E_2(u)=-a_1u$. Moreover,
\eqref{eq:codazzi} gives $E_2(\lambda)=-2\lambda a_1$ and $E_2(\mu)=-2\mu a_1$.
Differentiating the first relation of \eqref{eq:b} in the direction of $E_2$, we obtain
\[
 E_2(b_1)=\frac{\lambda-2}{\mu}\,a_1u .
\]
On the other hand, the Ricci equation \eqref{eq:ricci} gives
$E_2(b_1)=-a_1b_1=(2+\lambda)a_1u/\mu$. Comparing the two expressions, we obtain
$4a_1u/\mu=0$ and, since $u\neq0$, $a_1=0$.

Then $E_2(\lambda)=E_2(\mu)=E_2(u)=0$, so that $E_2(k_1)=E_2(k_2)=E_2(\alpha)=0$ and
$E_2(b_1)=0$, and, from \eqref{eq:gauss}, $E_2(K)=0$. Finally, from the first relation
of \eqref{eq:codazzi},
\[
 a_2=\frac{E_1(\lambda)+(1+2\lambda)u}{2\lambda},
\]
and, since $[E_2,E_1]=-a_2E_2$ and $E_2(\lambda)=0$, we have $E_2(E_1(\lambda))=0$. Thus
$E_2(a_2)=0$.
\end{proof}

\begin{lemma}\label{lem:coordinates}
Let $\varphi\colon M^2\to N^4(\eps)$ be a non-CMC surface with flat normal bundle such
that $\nabla^\perp E_3\neq0$ and $A_4\neq0$ at any point, and assume that \eqref{eq:340}
and \eqref{eq:342} hold for a frame $\{E_1,E_2\}$ as in Proposition~\ref{prop:frame}.
Then, around any point, there exist local coordinates $(s,t)$ and a positive function
$\rho=\rho(s)$ such that
\[
 E_1=\frac{\partial}{\partial s},\qquad E_2=\frac{1}{\rho(s)}\frac{\partial}{\partial t},
 \qquad g=ds^2+\rho(s)^2\,dt^2,\qquad \dot\rho=-a_2\rho ,
\]
where the dot denotes the derivative with respect to $s$, and the functions $f$, $k_1$,
$k_2$, $\alpha$, $a_2$, $b_1$ and $K$ depend only on $s$. Moreover, any function $F$ with
$E_2(F)=0$ depends only on $s$, $\partial/\partial t$ is a Killing vector field on $M$,
and the curves $t=\mathrm{const}$ are geodesics of $M$.
\end{lemma}

\begin{proof}
By \eqref{eq:340}, $\nabla_{E_1}E_1=-a_1E_2=0$ and $[E_1,E_2]=a_2E_2$. Let
$\{\omega^1,\omega^2\}$ be the dual coframe of $\{E_1,E_2\}$. Then
$d\omega^1(E_1,E_2)=-\omega^1([E_1,E_2])=0$, so, locally, $\omega^1=ds$ for some function
$s$, and $E_1(s)=1$, $E_2(s)=0$. Let $t$ be a function which is constant along the
integral curves of $E_1$ and whose restriction to an integral curve of $E_2$ is a
local diffeomorphism. Then $(s,t)$ are local coordinates, $E_1(s)=1$ and $E_1(t)=0$,
so $E_1=\partial/\partial s$, and, since $E_2(s)=0$, $E_2=g_{22}^{-1/2}\,\partial/\partial t$
and $g=ds^2+g_{22}\,dt^2$. If $E_2(F)=0$, then $\partial F/\partial t=0$, and $F$ depends
only on $s$ on a coordinate rectangle. By \eqref{eq:340} and \eqref{eq:342}, this applies
to $f$, $k_1$, $k_2$, $\alpha$, $a_2$, $b_1$ and $K$. Since
$[E_1,E_2]=-\left(\partial_s\log\sqrt{g_{22}}\right)E_2$, we get
$\partial_s\log\sqrt{g_{22}}=-a_2(s)$, so $\sqrt{g_{22}}=\psi(t)\,e^{-\int a_2\,ds}$ for a
positive function $\psi$. Replacing $t$ by a primitive of $\psi$, which does not change
the previous properties, we obtain the expressions of $E_2$ and $g$ in the statement,
with $\rho=e^{-\int a_2\,ds}$.
Finally, $\partial/\partial t$ is a Killing vector field, since the coefficients of $g$
do not depend on $t$, and the curves $t=\mathrm{const}$ are geodesics, since
$\nabla_{E_1}E_1=0$.
\end{proof}

\begin{remark}\label{rem:rotational}
Let $\varphi\colon M^2\to N^4(\eps)$ be a non-CMC surface with flat normal bundle, with
$\nabla^\perp E_3\neq0$ at any point and $A_4\neq0$ at any point, which satisfies
\eqref{eq:340} and \eqref{eq:342}, that is, a surface as in Lemma~\ref{lem:coordinates}
and Theorem~\ref{thm:41}. The geometric properties of these surfaces are given in
Proposition~\ref{prop:geometry} of Appendix~\ref{app:rotational}. In particular, for
$\eps=1$ they are rotational surfaces (Proposition~\ref{prop:geometry}~e)) and,
conversely, a curve in a totally geodesic $\Sph^3$ of $\Sph^4$ whose curvature and torsion
are those given in Proposition~\ref{prop:geometry}~e) determines a surface of this type
(Proposition~\ref{prop:converse2}). Moreover, for $\eps=1$, under some open conditions,
these surfaces are biharmonic if and only if the differential system \eqref{eq:Bsystem},
for the fourth component, the curvature and the torsion of the profile curve, is satisfied
(Proposition~\ref{prop:biharmonic-rot3}).
\end{remark}

By Lemma~\ref{lem:coordinates}, $f$ and $K$ are functionally dependent, so the surfaces
of Lemma~\ref{lem:coordinates} are W-surfaces; in particular, by
Proposition~\ref{prop:reduction}, on $M'$ the surfaces of
Proposition~\ref{prop:reduction} are W-surfaces. Moreover, \eqref{eq:340} and
\eqref{eq:342} are exactly (3.40)--(3.42) of \cite{AMOS}. We obtain the following version
of \cite[Theorem~4.1]{AMOS} without the W-hypothesis and the hypothesis
$3f^2+K-\eps\neq0$.

\begin{theorem}\label{thm:41}
Let $\varphi\colon M^2\to N^4(\eps)$ be a non-CMC surface with flat normal bundle. Assume
that $\nabla^\perp E_3\neq0$ and $A_4\neq0$ at any point, and that \eqref{eq:340} and
\eqref{eq:342} hold. Let $(s,t)$ be local coordinates as in Lemma~\ref{lem:coordinates}.
Then, $M$ is biharmonic if and only if the following ODE system is satisfied:
\begin{equation}\label{eq:41}
\left\{
\begin{aligned}
&\dot a_2=\eps+k_2\left(2f-k_2\right)-\alpha^2+a_2^2\\
&\left(3f-k_2\right)\dot f=-f\alpha b_1\\
&\dot\alpha=2\alpha a_2+k_2b_1\\
&\dot k_2=2a_2\left(k_2-f\right)-\alpha b_1\\
&\ddot f=a_2\dot f+f\left(b_1^2+k_1^2+k_2^2-2\eps\right)\\
&\dot b_1=-\frac{2b_1}{f}\dot f+a_2b_1-\alpha\left(k_2-k_1\right)
\end{aligned}\right. ,
\end{equation}
where the dot denotes the derivative with respect to $s$.
\end{theorem}

\begin{proof}
By Lemma~\ref{lem:coordinates}, $E_1=\partial/\partial s$ and the functions $f$, $k_1$,
$k_2$, $\alpha$, $a_2$, $b_1$ and $K$ depend only on $s$. First, the first, the third and
the fourth equations of \eqref{eq:41} hold for any such surface: replacing $a_1=0$ in
\eqref{eq:gauss}, we obtain the first one, and the third and the fourth ones are the third
and the first relations of \eqref{eq:codazziAMOS}, respectively.

By Proposition~\ref{prop:frame}, $M$ is biharmonic if and only if \eqref{eq:bicons} and
\eqref{eq:normalframe} hold. Since $E_2(f)=0$ and $b_2=0$, the second relation of
\eqref{eq:bicons} is identically satisfied, and the first one is the second equation of
\eqref{eq:41}. Replacing $a_1=0$, $b_2=0$, $E_2(f)=0$, $E_1=\partial/\partial s$ and the
fact that all the functions depend only on $s$ in \eqref{eq:normalframe}, we obtain the
last two equations of \eqref{eq:41}. Thus $M$ is biharmonic if and only if \eqref{eq:41}
holds.
\end{proof}

\begin{remark}
System \eqref{eq:41} is system (4.1) of \cite{AMOS}. Thus, by
Proposition~\ref{prop:reduction} and Theorem~\ref{thm:41}, (4.1), and its first order version (4.2), of \cite{AMOS}
characterize all non-CMC, non-PNMC biharmonic surfaces
with flat normal bundle in $N^4(\eps)$, not only the W-surfaces satisfying
$3f^2+K-\eps\neq0$.
\end{remark}

The hypothesis $\grad f\neq0$, contained in the non-CMC assumption, is not essential to
reach system \eqref{eq:41}. More precisely, from the hypotheses of
Theorem~\ref{thm:41} only the following ones are needed.

\begin{proposition}\label{prop:41minimal}
Let $\varphi\colon M^2\to N^4(\eps)$ be a surface with flat normal bundle such that
$H\neq0$ and $A_3\neq f\Id$ at any point, and let $\{E_1,E_2\}$ be a frame field as in
Proposition~\ref{prop:frame}. Assume that
\begin{equation}\label{eq:41hyp}
 a_1=0,\qquad b_2=0,\qquad E_2(f)=0 .
\end{equation}
Then, around any point, there exist local coordinates $(s,t)$ such that
$E_1=\partial/\partial s$ and the functions $f$, $k_1$, $k_2$, $\alpha$, $a_2$ and $b_1$
depend only on $s$. Moreover, $M$ is biharmonic if and only if system \eqref{eq:41} is
satisfied, where the dot denotes the derivative with respect to $s$.
\end{proposition}

\begin{proof}
The proof of Proposition~\ref{prop:frame} uses only that the normal bundle is flat,
$H\neq0$ and $A_3\neq f\Id$, so the frame $\{E_1,E_2\}$ exists and \eqref{eq:ricci},
\eqref{eq:gauss}, \eqref{eq:codazziAMOS}, \eqref{eq:bicons} and \eqref{eq:normalframe}
hold. By \eqref{eq:41hyp}, the second and the fourth relations of \eqref{eq:codazziAMOS}
give $E_2(k_1)=a_1(k_2-k_1)+\alpha b_2=0$ and $E_2(\alpha)=-2\alpha a_1-k_1b_2=0$, and
\eqref{eq:ricci} gives $E_2(b_1)=E_1(b_2)-a_1b_1-a_2b_2=0$. Then
$E_2(k_2)=2E_2(f)-E_2(k_1)=0$. Since $[E_1,E_2]=a_2E_2$, any function $F$ with
$E_2(F)=0$ satisfies $E_2(E_1(F))=E_1(E_2(F))-[E_1,E_2](F)=0$. Applying $E_2$ to the
first relation of \eqref{eq:codazziAMOS}, $E_1(k_2)=a_2(k_2-k_1)-\alpha b_1$, we obtain
$(k_2-k_1)E_2(a_2)=0$ and, since $k_1\neq k_2$, $E_2(a_2)=0$. Finally, from
\eqref{eq:gauss}, $E_2(K)=0$. Thus \eqref{eq:342} holds.

The proof of Lemma~\ref{lem:coordinates} uses only $a_1=0$, $E_2(f)=0$ and
\eqref{eq:342}, so it gives the local coordinates $(s,t)$ in the statement. The proof of
Theorem~\ref{thm:41} uses only these coordinates, $a_1=0$, $b_2=0$ and $E_2(f)=0$, so $M$
is biharmonic if and only if \eqref{eq:41} holds.
\end{proof}

The hypothesis $A_3\neq f\Id$ excludes the pseudo-umbilical surfaces. However,
for rotational surfaces, system \eqref{eq:41} still characterizes the biharmonicity in the
pseudo-umbilical case.

\begin{proposition}\label{prop:pseudoumbilical}
Let $\varphi\colon M^2\to\Sph^4$ be a rotational surface, that is, up to an isometry of
$\Sph^4$,
\[
 \varphi(s,t)=\left(\sigma_1(s),\sigma_2(s),\sigma_3(s),\rho(s)\cos t,\rho(s)\sin t\right),
\]
where $\sigma=(\sigma_1,\sigma_2,\sigma_3,\rho)$ is a curve parametrized by arc length in
$\Sph^3\cap\{x_4>0\}$. Assume that $\varphi$ is pseudo-umbilical, that is, $A_3=f\Id$, and
that $H\neq0$ at any point. Then $M$ has flat normal bundle and, with respect to the frame
$E_1=\partial/\partial s$, $E_2=\rho^{-1}\partial/\partial t$, the relations
\eqref{eq:frame} hold with $k_1=k_2=f$, and, with the notations \eqref{eq:connection} and
\eqref{eq:normalconnection}, \eqref{eq:41hyp} holds and the functions $f$, $\alpha$, $a_2$
and $b_1$ depend only on $s$. Moreover, $M$ is biharmonic if and only if system
\eqref{eq:41}, with $\eps=1$ and $k_2=f$, is satisfied, where the dot denotes the
derivative with respect to $s$.
\end{proposition}

\begin{proof}
Let $R_t$ be the rotation of angle $t$ in the plane $x_4x_5$, and let $J$ be the rotation of
angle $\pi/2$ in the same plane, so that $\varphi(s,t)=R_t\sigma(s)$, where we identify
$\sigma$ with $(\sigma_1,\sigma_2,\sigma_3,\rho,0)$. Then $\partial\varphi/\partial s=R_t\dot\sigma$ and
$\partial\varphi/\partial t=R_tJ\sigma=\rho R_t\partial_5$, so $E_1=R_t\dot\sigma$ and
$E_2=R_t\partial_5$. Let $\nu$ be a vector field along $\sigma$, with values in
$\R^4=\{x_5=0\}$, orthogonal to $\sigma$ and $\dot\sigma$. Then $\tilde\nu=R_t\nu$ is a normal
vector field along $\varphi$, the normal bundle is spanned by such vector fields, and
\[
 \frac{\partial\tilde\nu}{\partial t}=R_tJ\nu=\langle\nu,\partial_4\rangle R_t\partial_5
 =\langle\nu,\partial_4\rangle E_2 .
\]
Thus $\nabla^\perp_{E_2}\tilde\nu=0$ and $A_{\tilde\nu}E_2=-\rho^{-1}\langle\nu,\partial_4\rangle E_2$.
Hence all the shape operators are diagonal with respect to $\{E_1,E_2\}$, so they commute,
and, by \eqref{eq:generalricci}, the normal bundle is flat.

Since the rotations $R_t$ preserve $\varphi(M)$, they preserve $H$, so $E_3=H/|H|$ and $f=|H|$
are invariant, and we can choose $E_4$ invariant. Therefore $E_3=R_t\nu_3$ and $E_4=R_t\nu_4$,
with $\nu_3$, $\nu_4$ as above, and all the functions $f$, $\alpha$, $a_2$, $b_1$ depend only
on $s$; in particular, $E_2(f)=0$ and, by the above, $b_2=0$. Moreover, $A_3=f\Id$ and,
since $A_4$ is diagonal and $\tr A_4=0$, \eqref{eq:frame} holds with $k_1=k_2=f$. The
$E_2$-component of $\nabla^0_{E_1}E_1=R_t\ddot\sigma$ is $\langle\ddot\sigma,\partial_5\rangle=0$, so
$a_1=0$.

The proof of Proposition~\ref{prop:frame} uses only that $\{E_1,E_2\}$ diagonalizes $A_3$
and $A_4$ and that $\tr A_4=0$, so \eqref{eq:ricci}, \eqref{eq:gauss},
\eqref{eq:codazziAMOS}, \eqref{eq:bicons} and \eqref{eq:normalframe} hold with
$k_1=k_2=f$. By \eqref{eq:41hyp}, \eqref{eq:gauss} and the first and the third relations
of \eqref{eq:codazziAMOS} are the first, the fourth and the third equations of
\eqref{eq:41}, with $\eps=1$ and $k_2=f$, while the other relations of
\eqref{eq:codazziAMOS} and \eqref{eq:ricci} are identically satisfied. As in the proof of
Theorem~\ref{thm:41}, \eqref{eq:bicons} and \eqref{eq:normalframe} reduce to the second,
the fifth and the sixth equations of \eqref{eq:41}. Thus $M$ is biharmonic if and only if
\eqref{eq:41}, with $\eps=1$ and $k_2=f$, is satisfied.
\end{proof}

\subsection{On system \texorpdfstring{\eqref{eq:41}}{(4.3)}}\label{rem:compatible}
System \eqref{eq:41} alone, that is, without the hypotheses of Theorem~\ref{thm:41}, is
compatible. Indeed, let $\eps=1$, and let $A$, $c$ and $s_0$ be constants such that
$2A^2+c^2=1/8$. Then the functions
\begin{equation}\label{eq:41solutions}
 f=1,\qquad k_2=1,\qquad b_1=0,\qquad a_2=-\frac{\dot\rho}{\rho},\qquad
 \alpha=\frac{c}{\rho^2},
\end{equation}
where $\rho^2=\frac14+A\cos\left(2\sqrt2\left(s-s_0\right)\right)$, satisfy \eqref{eq:41}
on any open interval on which $\rho>0$. In fact, for $f=1$, $k_2=1$
and $b_1=0$ we have $k_1=1$, the second, the fourth, the fifth and the sixth equations of
\eqref{eq:41} are identically satisfied, and the first and the third equations reduce to
\[
 \dot a_2=2-\alpha^2+a_2^2,\qquad \dot\alpha=2\alpha a_2,
\]
whose general solution is given by \eqref{eq:41solutions}; moreover,
$\rho=\left(2+a_2^2+\alpha^2\right)^{-1/2}$. Since $k_1=k_2=f$ and $b_1=0$, these solutions
do not satisfy the hypotheses of Theorem~\ref{thm:41}. In Remark~\ref{rem:allsolutions}
we show that, on an open interval on which $f\neq0$, they are, up to the symmetry
$(f,k_2,\alpha)\mapsto(-f,-k_2,-\alpha)$ of \eqref{eq:41}, all the solutions of
\eqref{eq:41} with $\eps=1$.

Geometrically, these solutions, and so, by Remark~\ref{rem:allsolutions}, all the
solutions of \eqref{eq:41} with $\eps=1$ and $f\neq0$, represent minimal surfaces of a small
hypersphere $\Sph^3(2)$ of $\Sph^4$. Indeed, consider a solution \eqref{eq:41solutions}, the metric
$g=ds^2+\rho^2dt^2$, which satisfies $\dot\rho=-a_2\rho$, the frame
$E_1=\partial/\partial s$, $E_2=\rho^{-1}\partial/\partial t$, for which
$\nabla_{E_1}E_1=0$ and $[E_1,E_2]=a_2E_2$, and a normal bundle with an orthonormal frame
$\{E_3,E_4\}$ such that
\[
 A_3=\Id,\qquad A_4E_1=\alpha E_1,\qquad A_4E_2=-\alpha E_2,\qquad \nabla^\perp E_3=0 .
\]
Then the Ricci, Gauss and Codazzi equations, that is, \eqref{eq:ricci},
\eqref{eq:gauss} and \eqref{eq:codazziAMOS} with $a_1=0$, $b_1=b_2=0$ and $k_1=k_2=1$,
reduce to the first and the third equations of \eqref{eq:41}, so, by the fundamental
theorem of submanifolds, there exists, locally, a surface $\varphi\colon M^2\to\Sph^4$,
unique up to isometries of $\Sph^4$, with these data. Its normal bundle is flat and its
mean curvature vector field is $H=E_3$. Let $\nabla^0$ be the flat connection of $\R^5$.
For any vector field $X$ tangent to $M$, we have
$\nabla^0_X E_3=-A_3X+\nabla^\perp_XE_3=-X$, so $w=\varphi+E_3$ is a constant vector of
$\R^5$, with $\langle w,w\rangle=2$ and $\langle\varphi,w\rangle=1$. Thus $\varphi(M)$ lies
in $\Sph^4\cap\{x\in\R^5\mid\langle x,w\rangle=1\}$, which is a small hypersphere
$\Sph^3(2)$, of radius $1/\sqrt2$, and $E_3=w-\varphi$ is, along $\varphi$, the unit normal of
$\Sph^3(2)$ in $\Sph^4$. Since $H=E_3$, the mean curvature vector field of $\varphi$ in
$\Sph^3(2)$ vanishes, that is, $\varphi$ is a minimal surface of $\Sph^3(2)$. The solutions
with $f=-1$ give the same surfaces, with $E_3$ replaced by $-E_3$. We note that not every
minimal surface of $\Sph^3(2)$ is obtained in this way, but only those whose functions
$a_2$ and $\alpha$ depend, as above, only on one variable.

Moreover, these surfaces are rotational surfaces of $\Sph^4$. Indeed, the proofs of
Proposition~\ref{prop:geometry}~d) and e) use only $a_1=0$, $b_2=0$, the fact that all the
functions depend only on $s$, and \eqref{eq:Astructure} with $\eps=1$, which, in this case,
are the first, the fourth and the third equations of \eqref{eq:41}. Thus, up to an
isometry of $\Sph^4$, these surfaces are given by \eqref{eq:Aparam}, where
$\rho=\left(2+a_2^2+\alpha^2\right)^{-1/2}$ is the function in \eqref{eq:41solutions}.
The translation $s_0$ corresponds to a translation of the coordinate $s$, replacing $A$ by
$-A$ corresponds to replacing $s_0$ by $s_0+\pi/\left(2\sqrt2\right)$, and replacing $c$
by $-c$ corresponds to replacing $E_4$ by $-E_4$. Therefore, up to isometries of the domain
and of $\Sph^4$, the solutions \eqref{eq:41solutions} form a one-parameter family of
rotational surfaces in $\Sph^4$, parametrized by $|A|\in[0,1/4]$. Different values of
$|A|$ give non-congruent surfaces since, by \eqref{eq:gauss}, their Gaussian curvature is
$K=2-\alpha^2=2-c^2/\rho^4$, whose maximum, on the maximal interval on which $\rho>0$, is
$16|A|/(1+4|A|)$.

By Remark~\ref{rem:allsolutions}, it follows that no solution of \eqref{eq:41}, with
$\eps=1$, on an open interval on which $f\neq0$, satisfies all the hypotheses of
Proposition~\ref{prop:41minimal}. More precisely, the corresponding surfaces above have
flat normal bundle and $H\neq0$, and, with respect to the frame $\{E_1,E_2\}$, they
satisfy $a_1=0$, $b_2=0$ and $E_2(f)=0$, that is, \eqref{eq:41hyp}; the only hypothesis
of Proposition~\ref{prop:41minimal} which is not satisfied is $A_3\neq f\Id$, since
$k_1=k_2=f$, that is, $A_3=f\Id$ at any point.

For instance, for $c=0$ and $A=1/4$, we have $\alpha=0$ and
$\rho=\cos\left(\sqrt2\left(s-s_0\right)\right)/\sqrt2$, and we obtain the totally geodesic
sphere $\Sph^2(2)$ of $\Sph^3(2)$. For $c=\pm1/(2\sqrt2)$ and $A=0$, we have
$\rho=1/2$, $a_2=0$ and $\alpha=\pm\sqrt2$, and we obtain the Clifford torus of
$\Sph^3(2)$, the product of two circles of radius $1/2$, which is a minimal surface of
$\Sph^3(2)$.

\section{Proof of Theorem~\ref{thm:nonexistence}}\label{sec:conj}

In this section we prove Theorem~\ref{thm:nonexistence}. Let $\varphi\colon M^2\to\Sph^4$
be a biharmonic surface with flat normal bundle such that $H\neq0$, $\grad f\neq0$,
$A_3\neq f\Id$ and $\nabla^\perp E_3\neq0$ at any point. By
Proposition~\ref{prop:reduction} with $\eps=1$, the set $M'=\{p\in M\mid A_4(p)\neq0\}$ is
open and dense in $M$ and, around any point of $M'$, eventually changing the frame as in
Remark~\ref{rem:symmetry}, \eqref{eq:340}--\eqref{eq:342} hold. Thus, by
Lemma~\ref{lem:coordinates} and Theorem~\ref{thm:41}, there exist an open interval $I$ and
a solution of \eqref{eq:41} on $I$, with $\eps=1$, such that, by \eqref{eq:340},
\begin{equation}\label{eq:nondeg}
 f>0,\qquad \dot f\neq0,\qquad \alpha\neq0,\qquad b_1\neq0
\end{equation}
at any point of $I$. Here $\dot f=E_1(f)\neq0$, since $E_2(f)=0$ and $\grad f\neq0$.

\subsection{The normalized system}
We put
\begin{equation}\label{eq:newvars}
 \bar u=\frac{a_2}{f},\qquad \bar v=\frac{b_1}{f},\qquad \bar x=\frac{\alpha}{f},\qquad
 \bar y=\frac{k_2}{f},\qquad \bar z=\frac{\dot f}{f^2},\qquad e=\frac{1}{f^2},
\end{equation}
and we consider a new parameter $\tau$ defined by $d\tau/ds=f$. We denote by a prime the
derivative with respect to $\tau$.

\begin{lemma}\label{lem:normalized}
The functions \eqref{eq:newvars} satisfy the polynomial ODE system
\begin{equation}\label{eq:field}
\left\{
\begin{aligned}
&\bar u'=e+\bar y(2-\bar y)-\bar x^2+\bar u^2-\bar u\bar z\\
&\bar v'=\bar u\bar v+2\bar x(1-\bar y)-3\bar v\bar z\\
&\bar x'=2\bar x\bar u+\bar y\bar v-\bar x\bar z\\
&\bar y'=2\bar u(\bar y-1)-\bar x\bar v-\bar y\bar z\\
&\bar z'=\bar u\bar z+\bar v^2+(2-\bar y)^2+\bar y^2-2e-2\bar z^2\\
&e'=-2e\bar z
\end{aligned}\right.
\end{equation}
and the constraint
\begin{equation}\label{eq:C}
 C:=(3-\bar y)\,\bar z+\bar x\,\bar v=0 .
\end{equation}
Moreover, $\bar x\neq0$, $\bar v\neq0$, $\bar z\neq0$ and $\bar y\neq3$ at any point.
\end{lemma}

\begin{proof}
For a function $g=G/f^k$ we have $g'=f^{-1}\left(\dot G/f^k-kG\dot f/f^{k+1}\right)$.
Using \eqref{eq:41} with $\eps=1$ and $k_1=2f-k_2$, we obtain \eqref{eq:field}. For
example,
\[
 \bar u'=\frac{\dot a_2}{f^2}-\frac{a_2\dot f}{f^3}=e+\bar y(2-\bar y)-\bar x^2+\bar u^2-\bar u\bar z,
 \qquad
 \bar z'=\frac{\ddot f}{f^3}-\frac{2\dot f^2}{f^4}=\bar u\bar z+\bar v^2+(2-\bar y)^2+\bar y^2-2e-2\bar z^2,
\]
and $e'=-2\dot f/f^4=-2e\bar z$. The second equation of \eqref{eq:41} divided by $f^3$ is
\eqref{eq:C}. From \eqref{eq:nondeg} we get $\bar x\neq0$, $\bar v\neq0$ and
$\bar z\neq0$; if $\bar y=3$, then \eqref{eq:C} gives $\bar x\bar v=0$, which is a
contradiction.
\end{proof}

\begin{remark}\label{rem:equivalence}
The change of variables \eqref{eq:newvars}, together with $d\tau/ds=f$, gives an
equivalence between system \eqref{eq:41}, with $\eps=1$, and system \eqref{eq:field}
with the constraint \eqref{eq:C}, which does not use \eqref{eq:nondeg}. More precisely,
the proof of Lemma~\ref{lem:normalized} shows that any solution of \eqref{eq:41}, with
$\eps=1$ and $f>0$, gives a solution of \eqref{eq:field} and \eqref{eq:C} with $e>0$.
Conversely, let $(\bar u,\bar v,\bar x,\bar y,\bar z,e)$ be a solution of
\eqref{eq:field} and \eqref{eq:C}, defined on an open interval, with $e>0$. We put
$f=e^{-1/2}$, we consider the parameter $s$ defined by $ds/d\tau=1/f$, and we put
\[
 a_2=f\bar u,\qquad b_1=f\bar v,\qquad \alpha=f\bar x,\qquad k_2=f\bar y .
\]
Since $e'=-2e\bar z$, we have $f'=f\bar z$, so $\dot f=ff'=f^2\bar z$, where the dot
denotes the derivative with respect to $s$. Thus the functions \eqref{eq:newvars} are
recovered and, reversing the computations in the proof of Lemma~\ref{lem:normalized},
$(a_2,k_2,\alpha,f,b_1)$ is a solution of \eqref{eq:41}, with $\eps=1$ and $f>0$. By the
symmetry $(f,k_2,\alpha)\mapsto(-f,-k_2,-\alpha)$ of \eqref{eq:41}, the solutions with
$f<0$ are obtained in the same way. We note that the solutions of \eqref{eq:field} and
\eqref{eq:C} with $e\leq0$ do not correspond to solutions of \eqref{eq:41} with $\eps=1$.
\end{remark}

Let $D$ be the derivation of the polynomial ring
$\mathbb Q[\bar u,\bar v,\bar x,\bar y,\bar z,e]$ associated with the vector field
\eqref{eq:field}, that is, $D(F)$ is the derivative of $F$ along \eqref{eq:field}. Since
$C$ vanishes along the solution, we have $D^jC=0$ along the solution, for any $j\geq0$.

\subsection{The constraint polynomials}
A direct computation gives
\begin{multline*}
 DC=3\bar u\bar v\bar x-3\bar u\bar y\bar z+5\bar u\bar z+3\bar v^2-3\bar v\bar x\bar z-2\bar x^2\bar y+2\bar x^2-2\bar y^3+10\bar y^2\\
 +3\bar y\bar z^2-16\bar y-6\bar z^2+12+2e(\bar y-3).
\end{multline*}
We denote
\begin{equation}\label{eq:AB}
 \mathcal A(X,y)=X+(y-3)^2,\qquad
 \mathcal B(X,y,e)=(y-1)X+(y-3)\left(y^2-2y+2-e\right).
\end{equation}

\begin{lemma}\label{lem:uv}
Along the solution, we have
\begin{equation}\label{eq:uv}
 \bar v=-\frac{(3-\bar y)\bar z}{\bar x},\qquad
 \bar u=\frac{3\bar z^2\,\mathcal A\left(\bar x^2,\bar y\right)-2\bar x^2\,\mathcal B\left(\bar x^2,\bar y,e\right)}{4\bar x^2\bar z}.
\end{equation}
\end{lemma}

\begin{proof}
The first relation is \eqref{eq:C}. Replacing it in $DC$, we obtain an affine function of
$\bar u$, with the coefficient of $\bar u$ equal to $-4\bar z\neq0$. Solving $DC=0$ for
$\bar u$, we obtain the second relation.
\end{proof}

For $j=2,3,4,5$, we replace \eqref{eq:uv} in $D^jC$, take the numerator of the
resulting rational function, and remove the factors which are powers of $\bar x$ and
$\bar z$ and the integer content. The resulting polynomial is even in $\bar x$ and in
$\bar z$. Denoting $X=\bar x^2$, $Z=\bar z^2$ and $y=\bar y$, we obtain polynomials
$\mathcal P_j\in\mathbb Z[X,Z,y,e]$ with the following properties.
\begin{center}
\begin{tabular}{@{}ccc@{}}
\toprule
 & number of terms & degrees in $(X,Z,y,e)$\\
\midrule
$\mathcal P_2$ & 52 & $(4,2,6,2)$\\
$\mathcal P_3$ & 152 & $(6,3,9,3)$\\
$\mathcal P_4$ & 359 & $(8,4,12,4)$\\
$\mathcal P_5$ & 720 & $(10,5,15,5)$\\
\bottomrule
\end{tabular}
\end{center}
For example,
\begin{equation}\label{eq:P2}
 \mathcal P_2=-3\mathcal A^2Z^2-2X\mathcal M Z+4X^2\left((y-1)X+y^3-3y^2+3y+1-2e\right)\mathcal B,
\end{equation}
where
\begin{multline*}
 \mathcal M=(2y-3)X^2+\left(e(y+2)+4y^3-19y^2+15y-3\right)X\\
 +e\,(y-3)^2(y-9)+2y^5-16y^4+43y^3-33y^2-45y+81 .
\end{multline*}
Let
\[
 \mathcal S=\left\{(X,Z,y,e)\in\R^4\ \middle|\ X>0,\ Z>0,\ y\neq3,\
 \mathcal P_2=\mathcal P_3=\mathcal P_4=\mathcal P_5=0\right\}.
\]
By Lemmas~\ref{lem:normalized} and \ref{lem:uv}, the curve
$\tau\mapsto\left(\bar x(\tau)^2,\bar z(\tau)^2,\bar y(\tau),e(\tau)\right)$ lies in
$\mathcal S$.

\subsection{Elimination}

\begin{proposition}\label{prop:finite}
The image of $\mathcal S$ under the projection $(X,Z,y,e)\mapsto e$ is a finite set.
\end{proposition}

\begin{proof}
We note that $\mathcal A>0$ on $\mathcal S$. The resultants with respect to $Z$ are given
by
\[
\begin{gathered}
 \Res_Z(\mathcal P_2,\mathcal P_3)=48X^{7}\mathcal A^2\mathcal B\mathcal Q_8,\qquad
 \Res_Z(\mathcal P_2,\mathcal P_4)=36864X^{10}\mathcal A^2\mathcal B^2\mathcal Q_{10},\\
 \Res_Z(\mathcal P_2,\mathcal P_5)=12288X^{12}\mathcal A^2\mathcal B^2\mathcal Q_{14},\qquad
 \Res_Z(\mathcal P_3,\mathcal P_4)=-9216X^{15}\mathcal A^2\mathcal B^3\mathcal Q_{16},
\end{gathered}
\]
where $\mathcal Q_k\in\mathbb Z[X,y,e]$ is irreducible and has degree $k$ in $X$; the
polynomials $\mathcal Q_8$, $\mathcal Q_{10}$, $\mathcal Q_{14}$, $\mathcal Q_{16}$ have
$420$, $748$, $1832$, $2634$ terms, respectively. Since a common zero of two polynomials
is a zero of their resultant, and $X>0$, $\mathcal A>0$ on $\mathcal S$, any point of
$\mathcal S$ satisfies $\mathcal B=0$ or
$\mathcal Q_8=\mathcal Q_{10}=\mathcal Q_{14}=\mathcal Q_{16}=0$.

\textit{Case $\mathcal B=0$.} If $y=1$, then $\mathcal B=-2(1-e)$, so $e=1$. If $y\neq1$,
then
\[
 X=X_{\mathcal B}:=-\frac{(y-3)\left(y^2-2y+2-e\right)}{y-1},
\]
and $y^2-2y+2-e\neq0$, since $X\neq0$. The numerator of
$\mathcal P_2\left(X_{\mathcal B},Z,y,e\right)$ is $Z(y-3)^2(y-1)(cZ+d)$, where
\[
 c=-3(y-1)(2y-e-1)^2,\qquad d=2\left(y^2-2y+2-e\right)q,
\]
\[
 q=8y^4-10y^3e-22y^3+3y^2e^2+26y^2e+19y^2-8ye^2-18ye-6y+7e^2-2e+3 .
\]
Thus $cZ+d=0$. If $c=0$, that is, $e=2y-1$, then $d=-8(y-3)^2(y-1)^3\neq0$, which is a
contradiction. Hence $c\neq0$ and $Z=Z_{\mathcal B}:=-d/c$; since $Z\neq0$, we also have
$q\neq0$. Replacing $X=X_{\mathcal B}$ and $Z=Z_{\mathcal B}$ in $\mathcal P_3$ and
$\mathcal P_4$, the numerators factor as
\[
\begin{gathered}
 24\,q\,(2y-e-1)^2(y-3)^3(y-1)^3\left(y^2-2y+2-e\right)^3h_3,\\
 -1152\,q^2(2y-e-1)^2(y-3)^4(y-1)^4\left(y^2-2y+2-e\right)^4h_4,
\end{gathered}
\]
where $h_3,h_4\in\mathbb Z[y,e]$ are irreducible and coprime, with $45$ and $62$ terms.
Thus $h_3=h_4=0$, and this system has finitely many solutions $(y,e)$.

\textit{Case $\mathcal B\neq0$.} Then
$\mathcal Q_8=\mathcal Q_{10}=\mathcal Q_{14}=\mathcal Q_{16}=0$. Let
$S_k=\Res_X(\mathcal Q_8,\mathcal Q_k)\in\mathbb Z[y,e]$, $k=10,14,16$. We have
\[
 \gcd\left(S_{10},S_{14},S_{16}\right)=\mathrm{const}\cdot(e-1)(y-1)^4(y-3)^{24}.
\]
If three polynomials in two variables have infinitely many common zeros in $\mathbb C^2$,
then they have a non-constant common factor. Thus the common zeros of $S_{10}$, $S_{14}$,
$S_{16}$ outside the lines $e=1$, $y=1$ and $y=3$ form a finite set. The line $y=3$ is
excluded on $\mathcal S$, and it remains to consider the line $y=1$. Let
$\mathcal P^{(1)}_j(X,Z,e)=\mathcal P_j(X,Z,1,e)$. The resultants
$\Res_Z\left(\mathcal P^{(1)}_2,\mathcal P^{(1)}_k\right)$, $k=3,4,5$, factor as
$\mathrm{const}\cdot X^{m_k}(X+4)^2(e-1)^{n_k}G_k$, with $(m_k,n_k)=(7,3),(10,4),(12,5)$
and $G_k$ irreducible. Since $y=1$ and $\mathcal B\neq0$, we have $e\neq1$, so
$G_3=G_4=G_5=0$, and $e$ is a common root of $\Res_X(G_3,G_4)$ and $\Res_X(G_3,G_5)$.
The greatest common divisor of these two polynomials is $\mathrm{const}\cdot(e-1)^5$,
which is a contradiction. Thus there are no points of $\mathcal S$ with $y=1$ and
$\mathcal B\neq0$.
\end{proof}

\begin{proof}[Proof of Theorem~\ref{thm:nonexistence}]
Assume that such a surface exists. As shown above, we obtain an interval $J$ and a
solution of \eqref{eq:field}--\eqref{eq:C} on $J$ such that
$\left(\bar x(\tau)^2,\bar z(\tau)^2,\bar y(\tau),e(\tau)\right)\in\mathcal S$ for any
$\tau\in J$. By Proposition~\ref{prop:finite}, the continuous function $e$ takes values in
a finite set, and, since $J$ is connected, $e$ is constant. Since $e=1/f^2$, the function
$f$ is constant on an open subset, which contradicts $\dot f\neq0$.
\end{proof}

\begin{remark}
The argument of this section works for any $\eps\neq0$, with $e=\eps/f^2$, since $e$
constant implies that $f$ is constant. Together with Sections~\ref{sec:principal} and
\ref{sec:ODE}, this shows that Theorem~\ref{thm:nonexistence} holds in $N^4(\eps)$ for
any $\eps\neq0$. For $\eps\leq0$, the non-existence also follows from
\cite[Theorem~1.2]{FHY}.
\end{remark}

\begin{remark}\label{rem:allsolutions}
Using the argument of this section, we can
determine all the solutions of \eqref{eq:41}, with $\eps=1$, on an open interval on which
$f\neq0$, without assuming the hypotheses of Theorem~\ref{thm:41}. By the symmetry
$(f,k_2,\alpha)\mapsto(-f,-k_2,-\alpha)$ of \eqref{eq:41}, we can assume that $f>0$. The
first, the third, the fourth, the fifth and the sixth equations of \eqref{eq:41} give
$\dot a_2$, $\dot\alpha$, $\dot k_2$, $\ddot f$ and $\dot b_1$ as analytic functions of
$(a_2,\alpha,k_2,f,\dot f,b_1)$, so any solution is analytic, and each of the functions
$\dot f$, $\alpha$ and $b_1$ either vanishes identically or its zeros are isolated.

Assume that $\dot f$ does not vanish identically. If neither $\alpha$ nor $b_1$ vanishes
identically, then \eqref{eq:nondeg} holds on an open subinterval, and the argument of
Section~\ref{sec:conj} shows that $f$ is constant, which is a contradiction. If
$\alpha=0$, the second equation of \eqref{eq:41} gives $k_2=3f$, so $k_1=-f$, and the third
equation gives $b_1=0$. Then the fourth equation gives $a_2=3\dot f/(4f)$ and, denoting
$p=\dot f$, the fifth equation becomes
\[
 \dot p=\frac{3p^2}{4f}+10f^3-2f .
\]
Replacing $a_2$ and $\dot p$ in the first equation, which becomes
$\dot a_2=1-3f^2+a_2^2$, we obtain $p^2=14f^4-\frac{10}{3}f^2$. Differentiating this
relation and dividing by $2p$, on an open subinterval on which $p\neq0$, we get
$\dot p=28f^3-\frac{10}{3}f$, while the expression of $\dot p$ above gives
$\dot p=\frac{41}{2}f^3-\frac92f$. Comparing, we obtain $f^2=-7/45$, which is a
contradiction. If $b_1=0$ and $\alpha$ does
not vanish identically, the second equation gives again $k_2=3f$, and then the sixth
equation gives $4f\alpha=0$, which is a contradiction.

Therefore $f$ is constant, and the second equation of \eqref{eq:41} gives $\alpha b_1=0$.
If $b_1$ does not vanish identically, then $\alpha=0$, the third equation gives $k_2=0$,
the fourth equation gives $a_2=0$, and the first equation becomes $0=1$, which is a
contradiction. Thus $b_1=0$, and the fifth equation gives $k_1^2+k_2^2=2$. If $\alpha$ does
not vanish identically, the sixth equation gives $k_1=k_2$. If $\alpha=0$, then $k_1$ and
$k_2$ are constant, since $k_1+k_2=2f$ and $k_1k_2=2f^2-1$, and the fourth equation gives
$a_2(k_2-f)=0$; since $a_2=0$ would give, by the first equation, $2f^2=0$, we get $k_2=f$.
In both cases
$k_1=k_2=f$, so $f=1$, and the solution is one of the solutions \eqref{eq:41solutions} of
Subsection~\ref{rem:compatible}.

In particular, any solution of \eqref{eq:41}, with $\eps=1$ and $f\neq0$, has constant mean
curvature $f=\pm1$; the hypotheses of Theorem~\ref{thm:41} are used only to exclude these
solutions.
\end{remark}

\section{Proof of Theorem~\ref{thm:main}}\label{sec:proof}

As recalled in the introduction, the non-existence of PNMC biharmonic surfaces with
non-constant mean curvature follows from the classification of PNMC biconservative
surfaces in $\Sph^4$ obtained in \cite{NOTY}, and from the analogous results in $\R^4$ and
in $\mathbb H^4$ obtained in \cite{YT} and \cite{NR}, respectively. Since biharmonicity is
preserved by homotheties of the ambient space, we can state it in any space form
$N^4(\eps)$.

\begin{proposition}[\cite{NOTY,NR,YT}]\label{prop:PNMC}
There are no biharmonic surfaces in $N^4(\eps)$ such that, at any point, $H\neq0$,
$\grad f\neq0$, $A_3\neq f\Id$ and $\nabla^\perp E_3=0$.
\end{proposition}

\begin{proof}[Proof of Theorem~\ref{thm:main}]
Let $\varphi\colon M^2\to\Sph^4$ be a proper biharmonic surface with flat normal bundle.
We first show that $f$ is constant. Assume that $f^2=|H|^2$, which is a smooth function
on $M$, is not constant. Then there exists a non-empty open subset $U$ of $M$ on which
$d(f^2)\neq0$. Since $f^2\geq0$ attains its minimum at the zeros of $H$, we have $H\neq0$
on $U$, so that $f$ is smooth and $\grad f\neq0$ on $U$. By Lemma~\ref{lem:pseudo},
the open subset $U_1=\{p\in U\mid A_3(p)\neq f\Id\}$ is non-empty. If $\nabla^\perp E_3=0$
on a non-empty open subset of $U_1$, we obtain a contradiction with
Proposition~\ref{prop:PNMC}. Otherwise, the open subset of $U_1$ on which
$\nabla^\perp E_3\neq0$ is non-empty, and we obtain a contradiction with
Theorem~\ref{thm:nonexistence}. Therefore $f$ is constant and, since $\varphi$ is proper,
$f\neq0$. Thus $\varphi$ is a proper biharmonic CMC surface in $\Sph^4$ and, by the main
result of \cite{BO}, $\varphi(M)$ lies minimally in a small hypersphere $\Sph^3(2)$ of
$\Sph^4$.
\end{proof}

\begin{remark}
Conversely, a minimal surface of a small hypersphere $\Sph^3(2)$ is proper biharmonic in
$\Sph^4$ (see \cite{CMO}), and its normal bundle is flat, since the unit normal of
$\Sph^3(2)$ in $\Sph^4$, restricted to the surface, is parallel in the normal bundle. Thus
it satisfies all the hypotheses of Theorem~\ref{thm:main}.
\end{remark}

\appendix

\section{Rotational character of the surfaces in Theorem~\ref{thm:41}}\label{app:rotational}

In this appendix we consider the surfaces which satisfy the hypotheses of
Proposition~\ref{prop:41minimal}, without assuming that they are biharmonic. We describe their geometry and show that, for $\eps=1$, they
are, locally, rotational surfaces. Then, for $\eps=1$, we characterize the biharmonic ones
by means of their profile curves. We assume that $N^4(\eps)$ is simply connected.

For $\eps=1$, we consider $\Sph^4\subset\R^5$, and we denote by $(x_1,\ldots,x_5)$ the
canonical coordinates of $\R^5$, by $\{\partial_1,\ldots,\partial_5\}$ its canonical basis,
by $\nabla^0$ its flat connection, and by $R_t$ the rotation of angle $t$ in the plane
$x_4x_5$. We also consider $\Sph^3=\Sph^4\cap\{x_5=0\}\subset\R^4=\{x_5=0\}$; since it is
the intersection of $\Sph^4$ with a hyperplane through the origin, $\Sph^3$ is a great
hypersphere, so it is totally geodesic in $\Sph^4$. For a curve $\sigma$ parametrized by
arc length in $\Sph^3$, with curvature $\kappa_1>0$, $\{V_1,V_2,V_3\}$ denotes its Frenet
frame in $\Sph^3$, so that
\[
 \dot\sigma=V_1,\qquad \dot V_1=\kappa_1V_2-\sigma,\qquad \dot V_2=-\kappa_1V_1+\kappa_2V_3,
 \qquad \dot V_3=-\kappa_2V_2,
\]
where $\kappa_2$ is the torsion of $\sigma$ in $\Sph^3$ (the functions $\kappa_1$ and
$\kappa_2$ should not be confused with the principal curvatures $k_1$ and $k_2$).

\begin{proposition}\label{prop:geometry}
Let $\varphi\colon M^2\to N^4(\eps)$ be a surface with flat normal bundle such that
$H\neq0$ and $A_3\neq f\Id$ at any point, and let $\{E_1,E_2\}$ be a frame field as in
Proposition~\ref{prop:frame}. Assume that \eqref{eq:41hyp} holds. Then, locally, the
following hold.
\begin{enumerate}[label=\alph*)]
\item There exist local coordinates $(s,t)$ and a positive function $\rho=\rho(s)$ such
that
\begin{equation}\label{eq:Arho}
 E_1=\frac{\partial}{\partial s},\qquad E_2=\frac{1}{\rho(s)}\frac{\partial}{\partial t},
 \qquad g=ds^2+\rho(s)^2\,dt^2,\qquad \dot\rho=-a_2\rho,
\end{equation}
and the functions $f$, $k_1$, $k_2$, $\alpha$, $a_2$, $b_1$ and $K$ depend only on $s$. In
particular, $M$ is a W-surface, $\partial/\partial t$ is a Killing vector field on $M$, the
curves $t=\mathrm{const}$ are geodesics of $M$, and the coordinate curves are lines of
curvature of all the shape operators. Moreover, on the open subset where $\grad f\neq0$,
the curves $t=\mathrm{const}$ are integral curves of $\grad f/|\grad f|$ and the curves
$s=\mathrm{const}$ are the level curves of $f$.
\item The functions $f$, $k_2$, $\alpha$, $a_2$ and $b_1$ satisfy
\begin{equation}\label{eq:Astructure}
\left\{
\begin{aligned}
&\dot a_2=\eps+k_2\left(2f-k_2\right)-\alpha^2+a_2^2\\
&\dot k_2=2a_2\left(k_2-f\right)-\alpha b_1\\
&\dot\alpha=2\alpha a_2+k_2b_1
\end{aligned}\right. ,
\end{equation}
and $k_2\neq f$ at any point. Conversely, let $f$ and $b_1$ be functions on an open
interval $I$ such that $f>0$, and let $(a_2,k_2,\alpha)$ be a solution of
\eqref{eq:Astructure} on $I$ such that $k_2\neq f$. Then, locally, there exists a surface as in the
statement whose functions $f$, $k_2$, $\alpha$, $a_2$ and $b_1$ are the given ones, and it
is unique up to isometries of $N^4(\eps)$.
\item The curves $s=\mathrm{const}$ have constant curvature
$\varkappa=\sqrt{a_2^2+k_2^2+\alpha^2}$. If $\varkappa>0$, such a curve is a circle of
$N^4(\eps)$, that is, it also has zero torsion, and it lies in a totally geodesic surface
of $N^4(\eps)$; if $\varkappa=0$, it is a geodesic of $N^4(\eps)$.
\item There exists a one-parameter group $\{\Phi_c\}_{c\in\R}$ of isometries of
$N^4(\eps)$ such that $\varphi(s,t+c)=\Phi_c(\varphi(s,t))$.
\item Let $\eps=1$. Then the function $\rho^2\left(1+a_2^2+k_2^2+\alpha^2\right)$
is a positive constant $c^2$. Considering the coordinate $ct$ instead of $t$, and
$\rho/c$ instead of $\rho$ in \eqref{eq:Arho}, we have $c=1$ and, up to an isometry of
$\Sph^4$,
\begin{equation}\label{eq:Aparam}
 \varphi(s,t)=\bigl(\sigma_1(s),\sigma_2(s),\sigma_3(s),\rho(s)\cos t,\rho(s)\sin t\bigr),
\end{equation}
where $\sigma=(\sigma_1,\sigma_2,\sigma_3,\rho)$ is a curve parametrized by arc length in
$\Sph^3\cap\{x_4>0\}$, and $\Sph^3=\Sph^4\cap\{x_5=0\}$ is totally geodesic in $\Sph^4$. The
curvature of $\sigma$ in $\Sph^3$ is $\sqrt{k_1^2+\alpha^2}$ and, if $k_1^2+\alpha^2\neq0$ at
any point, the curvature $\kappa_1$ and the torsion $\kappa_2$ of $\sigma$ in $\Sph^3$ are
\begin{equation}\label{eq:Acurvatures}
 \kappa_1=\sqrt{k_1^2+\alpha^2}>0,\qquad
 \kappa_2=b_1+\frac{k_1\dot\alpha-\alpha\dot k_1}{k_1^2+\alpha^2}
 =\frac{2\left(fk_1b_1+f\alpha a_2-\alpha\dot f\right)}{k_1^2+\alpha^2},
\end{equation}
where the torsion is computed with respect to the Frenet frame $\{V_1,V_2,V_3\}$ of
$\sigma$ in $\Sph^3$ given, along $\sigma$, by $V_1=E_1$,
$V_2=\left(k_1E_3+\alpha E_4\right)/\kappa_1$ and $V_3=\left(-\alpha E_3+k_1E_4\right)/\kappa_1$.
Thus $\varphi(M)$ is a rotational surface: it
is obtained by rotating the profile curve $\sigma$ in the plane $x_4x_5$, about the
$3$-plane $x_4=x_5=0$. The profile curves are the curves $t=\mathrm{const}$, the
parallels are the curves $s=\mathrm{const}$, and $\rho$ is the Euclidean distance from the
$3$-plane $x_4=x_5=0$.
\end{enumerate}
\end{proposition}

\begin{proof}
a) By the proof of Proposition~\ref{prop:41minimal}, \eqref{eq:342} holds. The existence
of the coordinates, \eqref{eq:Arho}, the dependence on $s$ only of the functions, and the
facts that $\partial/\partial t$ is a Killing vector field and that the curves
$t=\mathrm{const}$ are geodesics are given by the proof of Lemma~\ref{lem:coordinates},
which uses only $a_1=0$, $E_2(f)=0$ and \eqref{eq:342}. Since $f$ and $K$ depend only on
$s$, we have $df\wedge dK=0$, so $M$ is a W-surface. The property of the coordinate curves
follows from \eqref{eq:frame}. Finally, on the open subset where $\grad f\neq0$, since
$E_2(f)=0$, we have $\dot f=E_1(f)\neq0$ and $E_1=\pm\grad f/|\grad f|$, and the curves
$s=\mathrm{const}$ are the level curves of $f$.

b) Replacing $a_1=0$ in \eqref{eq:gauss}, we obtain the first equation of
\eqref{eq:Astructure}. The second and the third equations are the first and the third
relations of \eqref{eq:codazziAMOS}, with $k_1=2f-k_2$. The other two relations of
\eqref{eq:codazziAMOS} and the Ricci equation \eqref{eq:ricci} are identically satisfied,
by \eqref{eq:41hyp} and \eqref{eq:342}. The condition $k_2\neq f$ means $k_1\neq k_2$.

Conversely, on $I\times\R$ we consider the metric $g$ given by \eqref{eq:Arho}, with
$\rho=e^{-\int a_2\,ds}$, the trivial vector bundle of rank two with an orthonormal frame
$\{E_3,E_4\}$, the symmetric bilinear form defined by
$B(E_1,E_1)=k_1E_3+\alpha E_4$, $B(E_1,E_2)=0$, $B(E_2,E_2)=k_2E_3-\alpha E_4$, with
$k_1=2f-k_2$, and the metric connection defined by $\nabla^\perp_{E_1}E_3=b_1E_4$,
$\nabla^\perp_{E_2}E_3=0$. The Gaussian curvature of $g$ is
$-\ddot\rho/\rho=\dot a_2-a_2^2$, and the curvature of $\nabla^\perp$ vanishes. Then the
Gauss, Codazzi and Ricci equations reduce to \eqref{eq:Astructure} and to identities, and
the fundamental theorem of submanifolds gives an immersion into $N^4(\eps)$, unique up to
isometries, with these data. Its mean curvature vector field is $H=fE_3\neq0$,
$A_3\neq f\Id$, since $k_1\neq k_2$, $\{E_1,E_2\}$ is a frame as in
Proposition~\ref{prop:frame}, and \eqref{eq:41hyp} holds.

c) Along a curve $s=\mathrm{const}$, let $T=E_2$. From the Gauss and Weingarten formulas,
\eqref{eq:connection}, \eqref{eq:normalconnection} and $a_1=b_2=0$, we get
\[
 \nabla^\varphi_{E_2}E_2=a_2E_1+k_2E_3-\alpha E_4,\quad
 \nabla^\varphi_{E_2}E_1=-a_2E_2,\quad
 \nabla^\varphi_{E_2}E_3=-k_2E_2,\quad
 \nabla^\varphi_{E_2}E_4=\alpha E_2 .
\]
Thus $\nabla^\varphi_TT=a_2E_1+k_2E_3-\alpha E_4$, whose norm $\varkappa$ satisfies
$E_2(\varkappa)=0$. If $\varkappa=0$, the curve is a geodesic of $N^4(\eps)$. If
$\varkappa>0$, let $\varkappa N=a_2E_1+k_2E_3-\alpha E_4$, with $|N|=1$. Then
$\nabla^\varphi_TT=\varkappa N$ and $\nabla^\varphi_TN=-\varkappa T$, with $\varkappa$
constant along the curve, so the curve is a circle lying in the totally geodesic surface
of $N^4(\eps)$ tangent to $\mathrm{span}\{T,N\}$.

d) For $c$ small, let $T_c(s,t)=(s,t+c)$, and consider the restrictions of $\varphi$ and
$\varphi\circ T_c$ to a connected coordinate rectangle. With respect to the frames
$\{E_1,E_2\}$ and $\{E_3,E_4\}$, these immersions have the same coefficients of the
metric, of the second fundamental form and of the normal connection, since these depend
only on $s$. By the uniqueness part of the fundamental theorem of submanifolds, there
exists an isometry $\Phi_c$ of $N^4(\eps)$ such that $\varphi\circ T_c=\Phi_c\circ\varphi$,
whose differential maps the frame $\{E_1,\ldots,E_4\}$ at $(s,t)$ into the frame at
$(s,t+c)$. Then $\Phi_{c+d}$ and $\Phi_c\circ\Phi_d$ coincide at a point, together with
their differentials, so they coincide. Thus $\{\Phi_c\}$ is a local one-parameter group of
isometries, generated by a Killing vector field $X$ of $N^4(\eps)$ with
$X\circ\varphi=\partial\varphi/\partial t$, and, since $N^4(\eps)$ is complete, it extends
to a one-parameter group.

e) Since $\langle\varphi,\varphi\rangle=1$,
we have $\nabla^0_XY=\nabla^\varphi_XY-\langle X,Y\rangle\varphi$ and $\nabla^0_X\varphi=X$.
Let
\[
 W=a_2E_1+k_2E_3-\alpha E_4-\varphi=\nabla^0_{E_2}E_2,\qquad
 w=|W|=\sqrt{1+a_2^2+k_2^2+\alpha^2}>0 .
\]
By \eqref{eq:connection} and $a_1=0$, $\nabla^0_{E_1}E_2=\nabla_{E_1}E_2+B(E_1,E_2)=0$.
Using
$\nabla^0_{E_1}E_1=k_1E_3+\alpha E_4-\varphi$, $\nabla^0_{E_1}E_3=-k_1E_1+b_1E_4$,
$\nabla^0_{E_1}E_4=-\alpha E_1-b_1E_3$ and \eqref{eq:Astructure} with $\eps=1$, a direct
computation gives
\[
\begin{aligned}
 \nabla^0_{E_1}W&=\left(\dot a_2-k_1k_2+\alpha^2-1\right)E_1
 +\left(a_2k_1+\dot k_2+\alpha b_1\right)E_3
 +\left(a_2\alpha+k_2b_1-\dot\alpha\right)E_4-a_2\varphi\\
 &=a_2W,
\end{aligned}
\]
and, by the computation in c), $\nabla^0_{E_2}W=-w^2E_2$. Thus the derivatives of $E_2$
and $W$ belong to $\mathrm{span}\{E_2,W\}$, and the $2$-plane
$\Pi=\mathrm{span}\{E_2,W\}$ is constant. Moreover, $\dot w=a_2w$ and $\dot\rho=-a_2\rho$,
so $\rho w$ is a positive constant $c$, and $\rho^2\left(1+a_2^2+k_2^2+\alpha^2\right)=c^2$.
After the change of coordinate and of $\rho$ in the statement, we have $\rho w=1$. Let
$\zeta=-W/w$. Then
\[
 \nabla^0_{\partial_s}\zeta=0,\qquad \nabla^0_{\partial_t}\zeta=\rho wE_2=E_2,\qquad
 \nabla^0_{\partial_t}E_2=\rho W=-\zeta,\qquad \nabla^0_{\partial_s}E_2=0,
\]
so $\{\zeta,E_2\}$ is an orthonormal basis of $\Pi$ which does not depend on $s$ and
rotates with unit angular speed in $t$. Choosing the coordinates of $\R^5$ such that $\Pi$
is the plane $x_4x_5$, and translating $t$, we have
$\zeta=(0,0,0,\cos t,\sin t)$ and $E_2=(0,0,0,-\sin t,\cos t)$. Since
$\langle\varphi,E_2\rangle=0$ and $\langle\varphi,W\rangle=-1$, the orthogonal projection
of $\varphi$ onto $\Pi$ is $\zeta/w=\rho(s)(0,0,0,\cos t,\sin t)$. Since
$\partial_t\varphi=\rho E_2\in\Pi$, the projection of $\varphi$ onto $\Pi^\perp$ depends
only on $s$. This gives \eqref{eq:Aparam}, with $\sigma(s)=\varphi(s,0)$, so
$|\sigma|=1$ and $|\dot\sigma|=|E_1|=1$.

Finally, at $t=0$ we have $E_2=\partial_5$, so $\sigma$, $E_1$, $E_3$ and $E_4$ are
orthogonal to $\partial_5$, that is, they lie in $\R^4=\{x_5=0\}$, and $\{E_1,E_3,E_4\}$ is,
along $\sigma$, an orthonormal frame of $T\Sph^3$. Along $\sigma$, we have $\dot\sigma=E_1$
and $\ddot\sigma=\nabla^0_{E_1}E_1=k_1E_3+\alpha E_4-\sigma$, so
$|\ddot\sigma+\sigma|=\sqrt{k_1^2+\alpha^2}$ is the curvature of $\sigma$ in $\Sph^3$. If
$k_1^2+\alpha^2\neq0$, then $\ddot\sigma+\sigma=\kappa_1V_2$, with $\kappa_1>0$ and $V_2$ as
in the statement. Moreover, by the expressions of $\nabla^0_{E_1}E_3$ and
$\nabla^0_{E_1}E_4$ above,
\[
 \frac{d}{ds}\left(k_1E_3+\alpha E_4\right)=-\kappa_1^2E_1+\left(\dot k_1-\alpha b_1\right)E_3
 +\left(\dot\alpha+k_1b_1\right)E_4 .
\]
Since $|V_2|=1$, it follows that $\dot V_2=-\kappa_1V_1+\kappa_2V_3$, with
\[
 \kappa_2=\frac{-\alpha\left(\dot k_1-\alpha b_1\right)+k_1\left(\dot\alpha+k_1b_1\right)}{\kappa_1^2}
 =b_1+\frac{k_1\dot\alpha-\alpha\dot k_1}{k_1^2+\alpha^2},
\]
and then $\dot V_3=-\kappa_2V_2$. The last expression of $\kappa_2$ in \eqref{eq:Acurvatures}
follows from \eqref{eq:Astructure}, with $k_1=2f-k_2$.
\end{proof}

We now show that the converse of Proposition~\ref{prop:geometry}~e) holds, when
$\alpha\neq0$. For $\eps=1$, Proposition~\ref{prop:converse2} below also gives an explicit
construction of the surface in the converse part of Proposition~\ref{prop:geometry}~b),
when $\alpha\neq0$. As in
\cite[Lemmas~2.1 and 2.2]{NR2}, we first relate the curvature and the torsion of a curve in
$\Sph^3$ to the existence of a constant vector which makes a prescribed inner product with
the position vector of the curve.

\begin{lemma}\label{lem:Aconst}
Let $\sigma\colon I\to\Sph^3\subset\R^4$ be a curve parametrized by arc length, with
curvature $\kappa_1>0$, torsion $\kappa_2$ and Frenet frame $\{V_1,V_2,V_3\}$. Assume that
there exists a constant vector $b\in\R^4$ such that
\begin{equation}\label{eq:Aconst}
 \langle\sigma(s),b\rangle=\Theta,\qquad \langle b,b\rangle=1,\qquad \langle b,V_3\rangle\neq0,
\end{equation}
where $\Theta\in C^\infty(I)$. Then
\begin{equation}\label{eq:Apositive}
 1-\Theta^2-\dot\Theta^2-\left(\frac{\ddot\Theta+\Theta}{\kappa_1}\right)^2>0
\end{equation}
and
\begin{equation}\label{eq:Atorsion}
 \kappa_2\langle b,V_3\rangle=\frac{d}{ds}\left(\frac{\ddot\Theta+\Theta}{\kappa_1}\right)
 +\kappa_1\dot\Theta .
\end{equation}
In particular, the torsion of $\sigma$ is uniquely determined by $\kappa_1$, $\Theta$ and the
sign of $\langle b,V_3\rangle$.
\end{lemma}

\begin{proof}
We write $b=\chi_0\sigma+\chi_1V_1+\chi_2V_2+\chi_3V_3$, where $\chi_i\in C^\infty(I)$,
$i\in\{0,\ldots,3\}$. Since $\langle\sigma,\sigma\rangle=1$, the relation
$\langle\sigma,b\rangle=\Theta$ gives $\chi_0=\Theta$. Differentiating
$\langle\sigma,b\rangle=\Theta$ and using $\dot\sigma=V_1$ and the fact that $b$ is
constant, we obtain $\chi_1=\langle b,V_1\rangle=\dot\Theta$. Differentiating again and
using $\dot V_1=\kappa_1V_2-\sigma$, we obtain $\ddot\Theta=\kappa_1\langle b,V_2\rangle-\Theta$,
so
\[
 \chi_2=\langle b,V_2\rangle=\frac{\ddot\Theta+\Theta}{\kappa_1}.
\]
Since $\langle b,b\rangle=1$ and $\chi_3=\langle b,V_3\rangle\neq0$, we get
$\chi_3^2=1-\Theta^2-\dot\Theta^2-\chi_2^2>0$, that is, \eqref{eq:Apositive}. Finally,
differentiating $\chi_2=\langle b,V_2\rangle$ and using $\dot V_2=-\kappa_1V_1+\kappa_2V_3$,
we obtain $\dot\chi_2=-\kappa_1\dot\Theta+\kappa_2\chi_3$, which is \eqref{eq:Atorsion}.
\end{proof}

A direct computation also gives the converse of Lemma~\ref{lem:Aconst}.

\begin{lemma}\label{lem:Aconstconverse}
Let $\sigma\colon I\to\Sph^3\subset\R^4$ be a curve parametrized by arc length, with
curvature $\kappa_1>0$, torsion $\kappa_2$ and Frenet frame $\{V_1,V_2,V_3\}$. Let
$\Theta\in C^\infty(I)$, and let $\ell\in C^\infty(I)$ be a function with $\ell\neq0$ at any
point such that
\begin{equation}\label{eq:Ahyp}
 \ell^2=1-\Theta^2-\dot\Theta^2-\left(\frac{\ddot\Theta+\Theta}{\kappa_1}\right)^2,
 \qquad
 \kappa_2\ell=\frac{d}{ds}\left(\frac{\ddot\Theta+\Theta}{\kappa_1}\right)+\kappa_1\dot\Theta .
\end{equation}
Then the vector
\begin{equation}\label{eq:Ab}
 b=\Theta\sigma+\dot\Theta V_1+\frac{\ddot\Theta+\Theta}{\kappa_1}V_2+\ell V_3
\end{equation}
is constant, and it is the unique constant vector $b\in\R^4$ which satisfies
\eqref{eq:Aconst} and $\langle b,V_3\rangle=\ell$.
\end{lemma}

\begin{proof}
Denote $m=(\ddot\Theta+\Theta)/\kappa_1$. Differentiating \eqref{eq:Ab} with the Frenet
equations, we obtain
\[
 \dot b=\left(\Theta+\ddot\Theta-\kappa_1m\right)V_1
 +\left(\kappa_1\dot\Theta+\dot m-\kappa_2\ell\right)V_2+\left(\kappa_2m+\dot\ell\right)V_3 .
\]
The coefficient of $V_1$ vanishes by the definition of $m$, and the coefficient of $V_2$
vanishes by the second relation of \eqref{eq:Ahyp}. Differentiating the first relation of
\eqref{eq:Ahyp}, we get $\ell\dot\ell=-\dot\Theta\left(\Theta+\ddot\Theta\right)-m\dot m
=-m\left(\kappa_1\dot\Theta+\dot m\right)=-\kappa_2m\ell$, so, since $\ell\neq0$, the
coefficient of $V_3$ also vanishes. Thus $b$ is constant, and the first relation of
\eqref{eq:Ahyp} gives $\langle b,b\rangle=1$. Clearly, $\langle\sigma,b\rangle=\Theta$ and
$\langle b,V_3\rangle=\ell$. The uniqueness follows from the proof of
Lemma~\ref{lem:Aconst}, which determines all the components of $b$ in the frame
$\{\sigma,V_1,V_2,V_3\}$.
\end{proof}

\begin{proposition}\label{prop:converse2}
Let $f$ and $b_1$ be functions on an open interval $I$ such that $f>0$, and let
$(a_2,k_2,\alpha)$ be a solution of \eqref{eq:Astructure}, with $\eps=1$, on $I$ such that $\alpha\neq0$ and $k_2\neq f$. Let $k_1=2f-k_2$,
\[
 \rho=\left(1+a_2^2+k_2^2+\alpha^2\right)^{-1/2},
\]
and let $\kappa_1>0$ and $\kappa_2$ be the functions given by \eqref{eq:Acurvatures}. Let
$\sigma\colon I\to\Sph^3$ be a curve parametrized by arc length with curvature $\kappa_1$
and torsion $\kappa_2$. Then, up to an isometry of $\Sph^3$,
$\sigma=(\sigma_1,\sigma_2,\sigma_3,\rho)$ is a curve in $\Sph^3\cap\{x_4>0\}$, and the map
$\varphi$ given by \eqref{eq:Aparam} satisfies all the hypotheses of
Proposition~\ref{prop:geometry}, with $\eps=1$. Moreover, with respect to the frames
$E_1=\partial/\partial s$, $E_2=\rho^{-1}\partial/\partial t$ and
\[
 E_3=\frac{1}{\kappa_1}R_t\left(k_1V_2-\alpha V_3\right),\qquad
 E_4=\frac{1}{\kappa_1}R_t\left(\alpha V_2+k_1V_3\right),
\]
where $\{V_1,V_2,V_3\}$ is the Frenet frame of $\sigma$ in $\Sph^3$, the functions $f$,
$k_2$, $\alpha$, $a_2$ and $b_1$ of $\varphi$ are the given ones, and $\varphi$ is unique up
to isometries of $\Sph^4$.
\end{proposition}

\begin{proof}
\textit{Step 1: the compatibility of the curvature and the torsion of $\sigma$ with the
condition that its fourth component is $\rho$.} First, by \eqref{eq:Astructure} with
$\eps=1$ and $k_1=2f-k_2$, a direct computation gives
\[
 a_2\dot a_2+k_2\dot k_2+\alpha\dot\alpha=a_2\left(1+a_2^2+k_2^2+\alpha^2\right),
\]
so $\dot\rho=-a_2\rho$. Then, using the first equation of \eqref{eq:Astructure},
$\ddot\rho=\left(a_2^2-\dot a_2\right)\rho=-\left(1+k_1k_2-\alpha^2\right)\rho$, that is,
\[
 \frac{\ddot\rho+\rho}{\kappa_1}=\frac{\left(\alpha^2-k_1k_2\right)\rho}{\kappa_1}.
\]
We apply Lemma~\ref{lem:Aconstconverse} with $\Theta=\rho$ and
\[
 \ell=\frac{2f\alpha\rho}{\kappa_1}\neq0 .
\]
Since $\rho^2\left(1+a_2^2+k_2^2+\alpha^2\right)=1$ and $\kappa_1^2=k_1^2+\alpha^2$, we have
\[
 1-\rho^2-\dot\rho^2-\left(\frac{\ddot\rho+\rho}{\kappa_1}\right)^2
 =\frac{\rho^2}{\kappa_1^2}\left(\left(k_2^2+\alpha^2\right)\left(k_1^2+\alpha^2\right)
 -\left(\alpha^2-k_1k_2\right)^2\right)
 =\frac{\rho^2\alpha^2\left(k_1+k_2\right)^2}{\kappa_1^2}=\ell^2,
\]
so the first relation of \eqref{eq:Ahyp} holds, and a direct computation, using
\eqref{eq:Astructure}, \eqref{eq:Acurvatures} and $\dot\rho=-a_2\rho$, gives
\[
 \frac{d}{ds}\left(\frac{\left(\alpha^2-k_1k_2\right)\rho}{\kappa_1}\right)-\kappa_1a_2\rho
 =\kappa_2\ell,
\]
that is, the second relation of \eqref{eq:Ahyp} holds. Thus there is no incompatibility
between the fact that $\sigma$ is determined, up to isometries of $\Sph^3$, by its
curvature $\kappa_1$ and its torsion $\kappa_2$, and the fact that its fourth component is
$\rho$. More precisely, by Lemma~\ref{lem:Aconstconverse},
\[
 b=\rho\left(\sigma-a_2V_1+\frac{\alpha^2-k_1k_2}{\kappa_1}V_2+\frac{2f\alpha}{\kappa_1}V_3\right)
\]
is the unique constant unit vector of $\R^4$ such that $\langle\sigma,b\rangle=\rho$ and
$\langle b,V_3\rangle=\ell$; in particular, $\langle\sigma,b\rangle>0$. Choosing the
coordinates of $\R^4$ such that $b=\partial_4$, we obtain that
$\sigma=(\sigma_1,\sigma_2,\sigma_3,\rho)$ is a curve in $\Sph^3\cap\{x_4>0\}$. Moreover, let
$\Psi$ be an isometry of $\Sph^3$, the restriction of an orthogonal transformation of
$\R^4$, and $\hat\sigma=\Psi\circ\sigma$. Then $\Psi$ maps the Frenet frame of $\sigma$ into
the Frenet frame of $\hat\sigma$ and, by the uniqueness in Lemma~\ref{lem:Aconstconverse},
it maps the vector $b$ of $\sigma$ into the vector $b$ of $\hat\sigma$. Since $\sigma$ is
unique up to isometries of $\Sph^3$, and these extend to isometries of $\Sph^4$ which fix
$\partial_5$, the map $\varphi$ is unique up to isometries of $\Sph^4$.

\textit{Step 2: $\varphi$ satisfies all the hypotheses of
Proposition~\ref{prop:geometry}, with $\eps=1$.} We have $\varphi(s,t)=R_t(\sigma(s))$, so
$\varphi_s=R_t(V_1)$ and $\varphi_t=\rho R_t(\partial_5)$. Thus $\varphi$ is an immersion
with $g=ds^2+\rho^2dt^2$, and $[E_1,E_2]=-(\dot\rho/\rho)E_2=a_2E_2$, so the function $a_2$
of $\varphi$ is the given one. Since $\nabla^0_{E_1}E_1=R_t(\ddot\sigma)=R_t(\kappa_1V_2-\sigma)$
is orthogonal to $\varphi_s$ and to $R_t(\partial_5)$, we get $\nabla_{E_1}E_1=0$, that is,
$a_1=0$. Moreover,
$\nabla^0_{E_2}E_1=\rho^{-1}\partial_tR_t(V_1)=\rho^{-1}\langle V_1,\partial_4\rangle R_t(\partial_5)
=(\dot\rho/\rho)E_2$ is tangent, so $B(E_1,E_2)=0$; hence all the shape operators are
diagonal in $\{E_1,E_2\}$, they commute, and, by \eqref{eq:generalricci}, the normal
bundle is flat. For any functions $x$, $y$ of $s$, the vector $R_t(xV_2+yV_3)$ is
orthogonal to $\varphi=R_t(\sigma)$ and to $\varphi_s=R_t(V_1)$, since $R_t$ is an isometry,
and to $\varphi_t$, since $R_t(xV_2+yV_3)$ belongs to the span of
$\partial_1,\partial_2,\partial_3$ and $R_t(\partial_4)$. Thus $\{E_3,E_4\}$ is an
orthonormal frame of the normal bundle. By the Frenet equations,
\[
 \nabla^\perp_{E_1}R_t(xV_2+yV_3)=R_t\left(\left(\dot x-\kappa_2y\right)V_2+\left(\dot y+\kappa_2x\right)V_3\right),
\]
and, since $\nabla^0_{E_2}R_t(xV_2+yV_3)=\rho^{-1}\langle xV_2+yV_3,\partial_4\rangle
R_t(\partial_5)$ is tangent, $\nabla^\perp_{E_2}R_t(xV_2+yV_3)=0$. Next, $B(E_1,E_1)$ is the
normal component of $\nabla^0_{E_1}E_1=R_t(\kappa_1V_2-\sigma)$, that is,
$B(E_1,E_1)=\kappa_1R_t(V_2)$, and $B(E_2,E_2)$ is the normal component of
$\nabla^0_{E_2}E_2=\rho^{-1}\partial_tR_t(\partial_5)=-\rho^{-1}R_t(\partial_4)$, that is,
$B(E_2,E_2)=-\rho^{-1}R_t\left(\langle\partial_4,V_2\rangle V_2+\langle\partial_4,V_3\rangle V_3\right)$,
where, by Step 1,
\[
 \langle\partial_4,V_2\rangle V_2+\langle\partial_4,V_3\rangle V_3
 =\frac{\rho}{\kappa_1}\left(\left(\alpha^2-k_1k_2\right)V_2+2f\alpha V_3\right).
\]
Since $R_t(V_2)=(k_1E_3+\alpha E_4)/\kappa_1$ and $R_t(V_3)=(-\alpha E_3+k_1E_4)/\kappa_1$, and
$k_1+k_2=2f$, we get
\[
 B(E_1,E_1)=k_1E_3+\alpha E_4,\qquad B(E_2,E_2)=k_2E_3-\alpha E_4 .
\]
Hence $H=fE_3$, so the mean curvature function of $\varphi$ is $f$, $E_3$ is the
normalized mean curvature vector field and $E_4$ is a unit normal vector field orthogonal
to $E_3$; moreover, \eqref{eq:frame} holds, with the given functions $k_1$, $k_2$ and
$\alpha$. Since $E_3$ and $E_4$ are of the form $R_t(xV_2+yV_3)$, we have $b_2=0$, and, by
the expression of $\nabla^\perp_{E_1}$ above, with $x=k_1/\kappa_1$ and
$y=-\alpha/\kappa_1$,
\[
 \langle\nabla^\perp_{E_1}E_3,E_4\rangle
 =\frac{\alpha\left(\dot x-\kappa_2y\right)+k_1\left(\dot y+\kappa_2x\right)}{\kappa_1}
 =\kappa_2-\frac{k_1\dot\alpha-\alpha\dot k_1}{k_1^2+\alpha^2}=b_1,
\]
by \eqref{eq:Acurvatures}. Therefore $H\neq0$; $E_2(f)=0$, since $f$ depends only on
$s$; and $A_3\neq f\Id$, since $k_1\neq k_2$. So $\varphi$ is a surface with flat normal
bundle such that $H\neq0$ and $A_3\neq f\Id$ at any point, $\{E_1,E_2\}$ is a frame as in
Proposition~\ref{prop:frame}, and \eqref{eq:41hyp} holds.
\end{proof}

Lemmas~\ref{lem:Aconst} and \ref{lem:Aconstconverse} and
Proposition~\ref{prop:converse2} give the converse of Proposition~\ref{prop:geometry}~e),
when $\alpha\neq0$.

The following proposition is an analogue of Proposition~\ref{prop:41minimal}: while
Proposition~\ref{prop:41minimal} characterizes the biharmonic surfaces by means of system
\eqref{eq:41} for the functions $f$, $k_2$, $\alpha$, $a_2$ and $b_1$, the biharmonic
surfaces \eqref{eq:Aparam} can be characterized directly by the fourth component, the
curvature and the torsion of the profile curve, as follows.

\begin{proposition}\label{prop:biharmonic-rot3}
Let $\sigma\colon I\to\Sph^3\cap\{x_4>0\}$ be a curve parametrized by arc length, with
fourth component $\rho$, curvature $\kappa_1>0$, torsion $\kappa_2$ and Frenet frame
$\{V_1,V_2,V_3\}$, and let $\varphi$ be given by \eqref{eq:Aparam}. Denote
\[
 m=\frac{\ddot\rho+\rho}{\kappa_1},\qquad \ell=\langle V_3,\partial_4\rangle,
\]
which satisfy
\begin{equation}\label{eq:Bcompat}
 \ell^2=1-\rho^2-\dot\rho^2-m^2,\qquad \kappa_2\ell=\dot m+\kappa_1\dot\rho .
\end{equation}
Assume that, at any point,
\begin{equation}\label{eq:Bopen}
 \ell\neq0,\qquad \kappa_1^2\rho^2\neq1-\rho^2-\dot\rho^2 .
\end{equation}
Then $\varphi$ satisfies all the hypotheses of Proposition~\ref{prop:geometry}, with
$\eps=1$, and $\varphi$ is biharmonic if and only if
\begin{equation}\label{eq:Bsystem}
\left\{
\begin{aligned}
&\rho^2\left(3\kappa_1\rho-m\right)\dot\kappa_1
 =\dot\rho\left(1-\rho^2-\dot\rho^2-3\kappa_1^2\rho^2-2\kappa_1\rho m\right)+\kappa_1\rho^2\kappa_2\ell\\
&\ddot\kappa_1=-\frac{2\dot\rho}{\rho}\,\dot\kappa_1
 +\kappa_1\left(\kappa_1^2+\kappa_2^2+\frac{\dot\rho^2}{\rho^2}-1\right)
 -\frac{m}{\rho^3}\left(2\kappa_1^2\rho^2+1-4\rho^2-2\dot\rho^2\right)\\
&\kappa_1\dot\kappa_2=-2\kappa_2\left(\dot\kappa_1+\frac{\kappa_1\dot\rho}{\rho}\right)
 +\frac{\ell}{\rho^3}\left(4\rho^2+2\dot\rho^2-1\right)
\end{aligned}\right. .
\end{equation}
Conversely, let $\rho>0$, $\kappa_1>0$, $\kappa_2$ and $\ell$ be functions on an open
interval $I$ such that, with $m=(\ddot\rho+\rho)/\kappa_1$, \eqref{eq:Bcompat},
\eqref{eq:Bopen} and \eqref{eq:Bsystem} hold. Then there exists a curve $\sigma$ in
$\Sph^3\cap\{x_4>0\}$, parametrized by arc length, with fourth component $\rho$, curvature
$\kappa_1$, torsion $\kappa_2$ and $\langle V_3,\partial_4\rangle=\ell$, and the
corresponding surface \eqref{eq:Aparam} is biharmonic; it is unique up to isometries of
$\Sph^4$.
\end{proposition}

\begin{proof}
By the proof of Lemma~\ref{lem:Aconst}, with $b=\partial_4$ and $\Theta=\rho$, which does
not use $\langle b,V_3\rangle\neq0$ except for the strict inequality \eqref{eq:Apositive},
we have $\langle\partial_4,V_1\rangle=\dot\rho$ and $\langle\partial_4,V_2\rangle=m$, and
\eqref{eq:Bcompat} holds. Assume now that \eqref{eq:Bopen} holds. On $I$ we consider the
functions
\begin{equation}\label{eq:Binverse}
\begin{gathered}
 f=\frac{\sqrt{\left(\kappa_1\rho-m\right)^2+\ell^2}}{2\rho},\qquad
 k_1=\frac{\kappa_1\left(\kappa_1\rho-m\right)}{2\rho f},\qquad k_2=2f-k_1,\qquad
 \alpha=\frac{\kappa_1\ell}{2\rho f},\\
 a_2=-\frac{\dot\rho}{\rho},\qquad
 b_1=\kappa_2-\frac{k_1\dot\alpha-\alpha\dot k_1}{\kappa_1^2} .
\end{gathered}
\end{equation}
Since $\ell\neq0$, we have $f>0$ and $\alpha\neq0$. Using \eqref{eq:Bcompat}, a direct
computation shows that $(a_2,k_2,\alpha)$ is a solution of \eqref{eq:Astructure}, with
$\eps=1$, that $\rho=\left(1+a_2^2+k_2^2+\alpha^2\right)^{-1/2}$,
$\kappa_1=\sqrt{k_1^2+\alpha^2}$ and $\ell=2f\alpha\rho/\kappa_1$, and that $\kappa_2$ is
given by \eqref{eq:Acurvatures}. Moreover,
\[
\begin{gathered}
 4\rho^2f^2=\kappa_1^2\rho^2-2\rho\ddot\rho-3\rho^2-\dot\rho^2+1,\qquad
 k_1-k_2=\frac{\kappa_1^2\rho^2-1+\rho^2+\dot\rho^2}{2\rho^2f},\\
 b_1=\frac{\kappa_1\kappa_2\rho\left(\kappa_1\rho-m\right)+\ell\left(\rho\dot\kappa_1+2\kappa_1\dot\rho\right)}
 {\left(\kappa_1\rho-m\right)^2+\ell^2},
\end{gathered}
\]
so \eqref{eq:Bopen} means $\alpha\neq0$ and $k_2\neq f$. Thus the functions $f$, $b_1$, $a_2$, $k_2$ and $\alpha$ satisfy the hypotheses of
Proposition~\ref{prop:converse2}.

Let $\tilde\sigma$ be the curve in $\Sph^3\cap\{x_4>0\}$, with fourth component $\rho$,
given by Proposition~\ref{prop:converse2} for these functions, and let $\tilde\varphi$ be
the corresponding surface \eqref{eq:Aparam}. The curves $\sigma$ and $\tilde\sigma$ have
the same curvature and torsion, so there exists an isometry $\Psi$ of $\Sph^3$, the
restriction of an orthogonal transformation of $\R^4$, such that $\sigma=\Psi\circ\tilde\sigma$,
which maps the Frenet frame of $\tilde\sigma$ into the Frenet frame of $\sigma$. By the
proof of Proposition~\ref{prop:converse2}, the constant unit vector $\Psi(\partial_4)$
satisfies $\langle\sigma,\Psi(\partial_4)\rangle=\langle\tilde\sigma,\partial_4\rangle=\rho$
and $\langle V_3,\Psi(\partial_4)\rangle=2f\alpha\rho/\kappa_1=\ell$. Since $\partial_4$
satisfies the same relations, by the uniqueness in Lemma~\ref{lem:Aconstconverse} we get
$\Psi(\partial_4)=\partial_4$. Thus $\Psi$ preserves $\mathrm{span}\{\partial_1,\partial_2,
\partial_3\}$, and, extending $\Psi$ to $\R^5$ by $\Psi(\partial_5)=\partial_5$, it commutes
with $R_t$, so $\varphi=\Psi\circ\tilde\varphi$. Therefore $\varphi$ satisfies all the
hypotheses of Proposition~\ref{prop:geometry}, with $\eps=1$, and its functions $f$,
$k_2$, $\alpha$, $a_2$ and $b_1$ are given by \eqref{eq:Binverse}.

By Proposition~\ref{prop:converse2}, $\tilde\varphi$ satisfies the hypotheses of
Proposition~\ref{prop:41minimal}, with the coordinates $(s,t)$, and, since the first, the
third and the fourth equations of \eqref{eq:41} are \eqref{eq:Astructure}, $\tilde\varphi$,
and so $\varphi$, is biharmonic if and only if the second, the fifth and the sixth equations
of \eqref{eq:41}, with $\eps=1$, hold. By means of \eqref{eq:Binverse} and
\eqref{eq:Bcompat}, a direct computation shows that the second equation of \eqref{eq:41}
is equivalent to the first equation of \eqref{eq:Bsystem}, and that, under this equation,
the fifth and the sixth equations of \eqref{eq:41} are equivalent to the last two equations
of \eqref{eq:Bsystem}.

Conversely, let $\sigma$ be a curve parametrized by arc length in $\Sph^3$ with curvature
$\kappa_1$ and torsion $\kappa_2$. By Lemma~\ref{lem:Aconstconverse}, with $\Theta=\rho$,
there exists a unique constant unit vector $b\in\R^4$ such that $\langle\sigma,b\rangle=\rho$
and $\langle b,V_3\rangle=\ell$. Choosing the coordinates of $\R^4$ such that $b=\partial_4$,
$\sigma$ is a curve in $\Sph^3\cap\{x_4>0\}$ with fourth component $\rho$ and
$\langle V_3,\partial_4\rangle=\ell$, and, by the first part, the surface \eqref{eq:Aparam} is
biharmonic. Since $\sigma$ is unique up to isometries of $\Sph^3$ and $b$ is unique, the
surface is unique up to isometries of $\Sph^4$.
\end{proof}

\begin{remark}
	Following this approach, it seems that it is not easier to prove that the solutions of
	system \eqref{eq:Bsystem} do not satisfy \eqref{eq:Bopen} than to prove that the solutions
	of \eqref{eq:41} do not satisfy the hypotheses of Proposition~\ref{prop:41minimal}, which
	is the key result of our paper.
\end{remark}

\begin{remark}\label{rem:examplesrot}
The examples at the end of Subsection~\ref{rem:compatible} can be described in the context of
Proposition~\ref{prop:biharmonic-rot3}. For the solutions \eqref{eq:41solutions} we have
$\rho^2\left(1+a_2^2+k_2^2+\alpha^2\right)=\rho^2\left(2+a_2^2+\alpha^2\right)=1$, so
$\rho$ is normalized as in item~e) of Proposition~\ref{prop:geometry}. We take $s_0=0$.

For $c=0$ and $A=1/4$, consider the curve
\[
 \sigma(s)=\frac{1}{\sqrt2}\left(1,\sin\left(\sqrt2s\right),0,\cos\left(\sqrt2s\right)\right),
 \qquad |s|<\frac{\pi}{2\sqrt2},
\]
in $\Sph^3\cap\{x_4>0\}$, parametrized by arc length, with fourth component
$\rho=\cos\left(\sqrt2s\right)/\sqrt2$. It is a circle of radius $1/\sqrt2$ of $\R^4$, and
its curvature and torsion in $\Sph^3$ are $\kappa_1=1$ and $\kappa_2=0$. Moreover,
$V_2=\ddot\sigma+\sigma$, $V_3=\pm\partial_3$, $m=-\rho$ and $\ell=0$. The corresponding
surface \eqref{eq:Aparam},
\[
 \varphi(s,t)=\frac{1}{\sqrt2}\left(1,\sin\left(\sqrt2s\right),0,
 \cos\left(\sqrt2s\right)\cos t,\cos\left(\sqrt2s\right)\sin t\right),
\]
is the totally geodesic sphere $\Sph^2(2)$ of the small hypersphere
$\Sph^4\cap\{x_1=1/\sqrt2\}$, which is a $\Sph^3(2)$.

For $c=\pm1/(2\sqrt2)$ and $A=0$, consider the curve
\[
 \sigma(s)=\left(\frac{1}{\sqrt2},\frac12\cos(2s),\frac12\sin(2s),\frac12\right)
\]
in $\Sph^3\cap\{x_4>0\}$, parametrized by arc length, with fourth component $\rho=1/2$. It
is a circle of radius $1/2$ of $\R^4$, and its curvature and torsion in $\Sph^3$ are
$\kappa_1=\sqrt3$ and $\kappa_2=0$. Moreover, $m=1/(2\sqrt3)$ and, choosing the orientation
of $V_3$ such that $\ell$ and $c$ have the same sign, $\ell=\pm\sqrt{2/3}$. The
corresponding surface \eqref{eq:Aparam},
\[
 \varphi(s,t)=\left(\frac{1}{\sqrt2},\frac12\cos(2s),\frac12\sin(2s),
 \frac12\cos t,\frac12\sin t\right),
\]
is the Clifford torus of the small hypersphere $\Sph^4\cap\{x_1=1/\sqrt2\}$.

In both cases \eqref{eq:Bcompat} holds, as it must, and a direct computation shows that
\eqref{eq:Bsystem} is also satisfied, in agreement with the fact that these surfaces are
minimal in $\Sph^3(2)$, and so biharmonic in $\Sph^4$ (see \cite{CMO}). Since
$\left(\kappa_1\rho-m,\ell\right)\neq(0,0)$, the formulas \eqref{eq:Binverse} still make
sense, and they give back $f=1$, $k_1=k_2=1$, $b_1=0$, $a_2=-\dot\rho/\rho$ and
$\alpha=c/\rho^2$, that is, the solutions \eqref{eq:41solutions}. However, \eqref{eq:Bopen}
does not hold: for both curves we have $\kappa_1^2\rho^2=1-\rho^2-\dot\rho^2$, which
corresponds to $k_2=f$, and, moreover, $\ell=0$ for the first curve, which corresponds to
$\alpha=0$. Thus, in agreement with Subsection~\ref{rem:compatible}, these surfaces do not
satisfy the hypotheses of Proposition~\ref{prop:geometry}, since $A_3=f\Id$, and
Proposition~\ref{prop:biharmonic-rot3} does not apply to them.
\end{remark}

\section{Mathematica code for the elimination}\label{app:code}

The following Wolfram Language code reproduces the algebraic elimination in
Proposition~\ref{prop:finite}. The symbols \texttt{x} and \texttt{z} stand
for $\bar x$ and $\bar z$, and \texttt{xx} and \texttt{zz} stand for
$X=\bar x^2$ and $Z=\bar z^2$. The program first computes $D^jC$ along
\eqref{eq:field} and then substitutes the expressions for $\bar u$ and
$\bar v$ from \eqref{eq:uv}. It removes integer content and the common
powers of $\bar x$ and $\bar z$, which are nonzero by \eqref{eq:nondeg}.
Polynomial gcds are compared up to nonzero rational constants.

For the last three resultants we use exact interpolation in $e$. If
$m=\deg_X F$ and $n=\deg_X G$, the Sylvester determinant gives the bound
$\deg_e\operatorname{Res}_X(F,G)\le n\deg_e F+m\deg_e G$.
The program uses one more evaluation point than this bound and checks that
both degrees in $X$ are preserved at every point. Evaluating the listing
in a fresh Mathematica kernel returns \texttt{True} when all checks pass.

%[basicstyle=\ttfamily\footnotesize,columns=fullflexible,keepspaces=true,showstringspaces=false,breaklines=true,breakatwhitespace=true]

\begin{lstlisting}[style=mma]
	Begin["AMOSElimination`"];
	ClearAll["AMOSElimination`*"];
	
	check[condition_, label_] := If[!TrueQ[condition],
	Print["FAILED: ", label]; Throw[$Failed, "AMOS"]];
	termCount[p_] := Length[CoefficientRules[p, Variables[p]]];
	primitive[p_] := Last[FactorTermsList[Expand[p]]];
	associates[p_, q_] := With[{ratio = Cancel[p/q]},
	IntegerQ[ratio] && ratio != 0 || Head[ratio] === Rational];
	quotient[p_, q_] := Module[{r = Cancel[p/q]},
	check[PolynomialQ[r, Variables[p]], "exact polynomial division"]; Expand[r]];
	lie[p_] := Expand[Sum[field[[i]] D[p, vars[[i]]], {i, Length[vars]}]];
	squaredConstraint[p_] := Module[{num, rules, powers, shift},
	num = Numerator[Together[p /. uvRules]];
	rules = CoefficientRules[num, {x,z,y,e}];
	powers = First /@ rules;
	shift = {Min[powers[[All,1]]], Min[powers[[All,2]]]};
	primitive[FromCoefficientRules[Map[Function[rule,
	With[{indices = (Take[First[rule],2]-shift)/2},
	check[VectorQ[indices,IntegerQ], "parity"];
	Join[indices,Drop[First[rule],2]] -> Last[rule]]], rules], {xx,zz,y,e}]]];
	
	run[] := Catch[
	vars = {u,v,x,y,z,e};
	field = {e+y (2-y)-x^2+u^2-u z, u v+2 x (1-y)-3 v z,
		2 x u+y v-x z, 2 u (y-1)-x v-y z,
		u z+v^2+(2-y)^2+y^2-2 e-2 z^2, -2 e z};
	aa = xx+(y-3)^2;
	bb = (y-1) xx+(y-3) (y^2-2 y+2-e);
	jets = NestList[lie,(3-y) z+x v,5];
	uvRules = {v -> (y-3) z/x, u -> ((3 z^2 aa-2 xx bb)/(4 xx z) /. xx->x^2)};
	check[(Together /@ (Take[jets,2] /. uvRules)) === {0,0}, "C and DC"];
	Do[pol[j] = squaredConstraint[jets[[j+1]]], {j,2,5}];
	check[Table[termCount[pol[j]],{j,2,5}] === {52,152,359,720}, "P support"];
	check[And@@Table[Exponent[pol[j],{xx,zz,y,e}] === {2 j,j,3 j,j},{j,2,5}], "P degrees"];
	mm = (2 y-3) xx^2 + (e (y+2)+4 y^3-19 y^2+15 y-3) xx +
	e (y-3)^2 (y-9)+2 y^5-16 y^4+43 y^3-33 y^2-45 y+81;
	check[associates[pol[2],-3 aa^2 zz^2-2 xx mm zz+
	4 xx^2 ((y-1) xx+y^3-3 y^2+3 y+1-2 e) bb], "displayed P2"];
	Print["P2 through P5 verified."];
	
	(* B=0. On y=1 one gets e=1. Otherwise X can be eliminated. *)
	check[Expand[(bb /. y->1)-2 (e-1)] === 0, "B on y=1"];
	t = y^2-2 y+2-e;
	q = 8 y^4-10 y^3 e-22 y^3+3 y^2 e^2+26 y^2 e+19 y^2-
	8 y e^2-18 y e-6 y+7 e^2-2 e+3;
	c = -3 (y-1) (2 y-e-1)^2; d = 2 t q;
	xb = -(y-3) t/(y-1);
	check[Expand[(d /. e->2 y-1)+8 (y-3)^2 (y-1)^3] === 0, "c=0"];
	check[associates[(y-1)^4 (pol[2] /. xx->xb), zz (y-3)^2 (y-1) (c zz+d)], "B=0 linear equation"];
	Do[
	factors = Rest[FactorList[Numerator[Together[pol[j] /. {xx->xb,zz->-d/c}]]]];
	residual = Select[factors, Function[factor,
	!AnyTrue[{q,2 y-e-1,y-3,y-1,t},Function[known,associates[First[factor],known]]]]];
	check[Length[residual] === 1 && residual[[1,2]] === 1, "B=0 residual"];
	h[j] = residual[[1,1]], {j,3,4}];
	check[{termCount[h[3]],termCount[h[4]]} === {45,62} &&
	associates[PolynomialGCD[h[3],h[4]],1], "coprime h3,h4"];
	
	(* B != 0 on the line y=1. *)
	Do[{j,m,n} = row;
	gg[j] = primitive[quotient[Resultant[pol[2] /. y->1,pol[j] /. y->1,zz],
	xx^m (xx+4)^2 (e-1)^n]], {row,{{3,7,3},{4,10,4},{5,12,5}}}];
	check[And @@ Table[IrreduciblePolynomialQ[gg[j]],{j,3,5}],
	"G irreducibility"];
	special = PolynomialGCD[Resultant[gg[3],gg[4],xx],Resultant[gg[3],gg[5],xx]];
	check[associates[special,(e-1)^5],"y=1 gcd"];
	Print["Exceptional branches verified."];
	
	Do[{i,j,k,m,n,content,count} = row;
	raw = quotient[Resultant[pol[i],pol[j],zz],xx^m aa^2 bb^n];
	check[Abs[First[FactorTermsList[raw]]] === content,"resultant content"];
	qq[k] = primitive[raw];
	check[Exponent[qq[k],xx] === k && termCount[qq[k]] === count,"Q support"];
	check[IrreduciblePolynomialQ[qq[k]],"Q irreducibility"];
	Print["Q",k," verified."],
	{row,{{2,3,8,7,1,48,420},{2,4,10,10,2,36864,748},
			{2,5,14,12,2,12288,1832},{3,4,16,15,3,9216,2634}}}];
	
	Do[ss[k] = exactResultant[qq[8],qq[k]];
	Print["S",k," degrees in (y,e): ",Exponent[ss[k],{y,e}]], {k,{10,14,16}}];
	common = PolynomialGCD[ss[10],ss[14],ss[16]];
	check[associates[common,(e-1) (y-1)^4 (y-3)^24],"generic gcd"];
	Print["Generic gcd up to a constant: ",Factor[primitive[common]]];
	True, "AMOS"];
	
	(* Evaluate the Sylvester determinant at enough integer e-values to recover
	the bivariate resultant. Its e-degree is at most n deg_e F + m deg_e G. *)
	exactResultant[f_,g_] := Module[{m=Exponent[f,xx],n=Exponent[g,xx],bound,
		nodes,values,fs,gs,differences,coeffs,result},
	bound = n Exponent[f,e]+m Exponent[g,e];
	nodes = -2-Range[0,bound];
	values = Table[fs = f /. e->node; gs = g /. e->node;
	check[Exponent[fs,xx] === m && Exponent[gs,xx] === n,"specialization degrees"];
	CoefficientList[Resultant[fs,gs,xx],y], {node,nodes}];
	values = PadRight[values];
	differences = values;
	coeffs = Table[With[{first=First[differences]},
	differences=Differences[differences]; first/((-1)^j j!)],{j,0,bound}];
	(* Convert Newton coefficients to ordinary e-powers for all y-coefficients. *)
	result = Fold[(Expand[#2[[1]]+(e-#2[[2]]) #1])&,Last[coeffs],
	Reverse[Transpose[{Most[coeffs],Most[nodes]}]]];
	FromDigits[Reverse[result],y] // Expand];
	
	End[];
	AMOSElimination`run[]
\end{lstlisting}


\begin{thebibliography}{99}

\bibitem{AMOS}
\c S.~Andronic, S.~Montaldo, C.~Oniciuc and A.~Sanna,
{Biconservative Weingarten surfaces with flat normal bundle in $N^4(\epsilon)$},
\emph{Math. Nachr.} \textbf{299} (2026), 1624--1664.

\bibitem{BMO}
A.~Balmu\c{s}, S.~Montaldo and C.~Oniciuc,
{Classification results for biharmonic submanifolds in spheres},
\emph{Israel J. Math.} \textbf{168} (2008), 201--220.

\bibitem{BalmusMontaldoOniciuc2013}
A.~Balmu\c{s}, S.~Montaldo and C.~Oniciuc,
{Biharmonic PNMC submanifolds in spheres},
\emph{Ark. Mat.} \textbf{51} (2013), no. 2, 197--221.

\bibitem{BO}
A.~Balmu\c{s} and C.~Oniciuc,
{Biharmonic surfaces of $\Sph^4$},
\emph{Kyushu J. Math.} \textbf{63} (2009), 339--345.

\bibitem{CMO}
R.~Caddeo, S.~Montaldo and C.~Oniciuc,
{Biharmonic submanifolds in spheres},
\emph{Israel J. Math.} \textbf{130} (2002), 109--123.

\bibitem{CaddeoMontaldoOniciucPiu2014}
R.~Caddeo, S.~Montaldo, C.~Oniciuc, and P.~Piu,
Surfaces in three-dimensional space forms with divergence-free stress-bienergy tensor,
\emph{Ann. Mat. Pura Appl. (4)}, \textbf{193}(2) (2014), 529--550.

\bibitem{Chen}
B.-Y.~Chen,
{Some open problems and conjectures on submanifolds of finite type},
\emph{Soochow J. Math.} \textbf{17} (1991), 169--188.

\bibitem{Chen2019}
B.-Y.~Chen,
{Chen’s Biharmonic Conjecture and Submanifolds with Parallel Normalized Mean Curvature Vector},
\emph{Mathematics} \textbf{7} (2019), no. 8, 710.

\bibitem{Dimitric}
I.~Dimitri\'c,
{Submanifolds of $E^m$ with harmonic mean curvature vector},
\emph{Bull. Inst. Math. Acad. Sinica} \textbf{20} (1992), 53--65.

\bibitem{FetcuOniciuc2022}
D.~Fetcu and C.~Oniciuc,
{Biharmonic and biconservative hypersurfaces in space forms}, 
In {\em Differential geometry and global analysis---in honor of {T}adashi {N}agano}, {\em Contemp. Math.} \textbf{777}, Amer. Math. Soc., [Providence], RI, (2022), 65--90. 

\bibitem{FHY}
Y.~Fu, M.-C.~Hong and D.~Yang,
{Biharmonic submanifolds in nonpositive space forms},
preprint 2026, arXiv:2611.04159v1.

\bibitem{Jiang}
G.~Y.~Jiang,
{$2$-harmonic maps and their first and second variational formulas},
\emph{Chinese Ann. Math. Ser. A} \textbf{7} (1986), 389--402.

\bibitem{Jiang2}
G.~Y.~Jiang,
{The conservation law for $2$-harmonic maps between Riemannian manifolds},
\emph{Acta Math. Sinica} \textbf{30} (1987), 220--225.

\bibitem{Lawson}
H.~B.~Lawson,
{Complete minimal surfaces in $S^3$},
\emph{Ann. of Math. (2)} \textbf{92} (1970), 335--374.

\bibitem{LMO}
E.~Loubeau, S.~Montaldo and C.~Oniciuc,
{The stress-energy tensor for biharmonic maps},
\emph{Math. Z.} \textbf{259} (2008), 503--524.

%\bibitem{MO}
%S.~Montaldo and C.~Oniciuc,
%{A short survey on biharmonic maps between Riemannian manifolds},
%\emph{Rev. Un. Mat. Argentina} \textbf{47} (2006), 1--22.

\bibitem{Maeta2026}
S.~Maeta, 
{Biharmonic rotational surfaces in the four-dimensional Euclidean space are minimal},
\emph{arXiv:2605.09587}

\bibitem{NOTY}
S.~Nistor, C.~Oniciuc, N.~C.~Turgay and R.~Ye\u{g}in \c{S}en,
{Biconservative surfaces in the $4$-dimensional Euclidean sphere},
\emph{Ann. Mat. Pura Appl. (4)} \textbf{202} (2023), 2345--2377.

\bibitem{NR}
S.~Nistor and M.~Rusu,
{Intrinsic characterizations of biconservative surfaces in the four-dimensional
hyperbolic space},
\emph{Mediterr. J. Math.} \textbf{21} (2024), Paper No.~225.

\bibitem{NR2}
S.~Nistor and M.~Rusu,
{Extrinsic characterizations of biconservative surfaces in the $4$-dimensional
hyperbolic space},
preprint 2026, arXiv:2605.22396v2.

\bibitem{Oniciuc}
C.~Oniciuc,
{Biharmonic maps between Riemannian manifolds},
\emph{An. \c Stiin\c t. Univ. Al. I. Cuza Ia\c si. Mat. (N.S.)}
\textbf{48} (2002), 237--248.

\bibitem{ChenOu2020Book}
Y.-L. Ou and B.-Y. Chen,
{Biharmonic submanifolds and biharmonic maps in {R}iemannian geometry},
\emph{World Scientific Publishing Co. Pte. Ltd.}, Hackensack, NJ, 2020.

\bibitem{YT}
R.~Ye\u{g}in \c{S}en and N.~C.~Turgay,
{On biconservative surfaces in $4$-dimensional Euclidean space},
\emph{J. Math. Anal. Appl.} \textbf{460} (2018), 565--581.

\end{thebibliography}
\end{document}